\documentclass[11pt]{amsart}
\usepackage[backref=page]{hyperref}
\usepackage{hyperref}
\usepackage{amssymb,amsfonts,textcomp}

\usepackage{enumerate}

\usepackage[utf8]{inputenc}
\usepackage[T1]{fontenc}

\usepackage[mathscr]{eucal}

\usepackage[a4paper, twoside=false, vmargin={2cm,3cm}, includehead]{geometry}

\newtheorem{lemma}{Lemma}[section]
\newtheorem{theorem}[lemma]{Theorem}
\newtheorem*{theorem*}{Theorem}

\newtheorem*{proposition*}{Proposition}

\newtheorem{conjecture}{Conjecture}
\newtheorem*{conjecture*}{Conjecture}
\newtheorem{problem}{Problem}
\newtheorem*{problem*}{Problem}

\theoremstyle{definition}
\newtheorem*{claim*}{Claim}

\newtheorem{definition}[lemma]{Definition}

\newcommand{\legendre}[2]{\ensuremath{\left( \frac{#1}{#2} \right) }}

\newcommand{\C}{{\mathbb C}}
\newcommand{\E}{{\mathbb E}}
\newcommand{\D}{{\mathbb D}}

\newcommand{\N}{{\mathbb N}}
\renewcommand{\P}{{\mathbb P}}
\newcommand{\Q}{{\mathbb Q}}
\newcommand{\R}{{\mathbb R}}
\renewcommand{\S}{\mathbb{S}}
\newcommand{\T}{{\mathbb T}}
\newcommand{\Z}{{\mathbb Z}}
\newcommand{\U}{{\mathbb U}}

\newcommand{\CA}{{\mathcal A}}

\newcommand{\CM}{{\mathcal M}}

\newcommand{\CP}{{\mathcal P}}

\newcommand{\CX}{{\mathcal X}}

\newcommand{\wh}{\widehat}

\newcommand{\norm}[1]{\left\Vert #1\right\Vert}

\newcommand{\inv}{^{-1}}

\newcommand{\er}{{\text{\rm er}}}

\DeclareMathOperator{\Imag}{Im}
\renewcommand{\Im}{\Imag}

\begin{document}

			\title[Partition regularity of homogeneous quadratics]{Partition regularity of homogeneous quadratics: Current trends and challenges}

	\author{Nikos Frantzikinakis}
	\address[Nikos Frantzikinakis]{University of Crete, Department of Mathematics and Applied Mathematics, Heraklion, Greece}
	\email{frantzikinakis@gmail.com}

	\begin{abstract}
		Suppose   we  partition the integers into finitely many cells. Can we always find a solution of
		the equation $x^2+y^2=z^2$ with  $x,y,z$ on the same cell? What about
		more general homogeneous quadratic equations in three variables? These are basic questions in  arithmetic Ramsey theory, which have recently been partially answered using ideas inspired by ergodic theory and tools such as
		Gowers-uniformity properties and concentration estimates of bounded multiplicative functions.		
		The aim of this article 		is to provide an introduction to this exciting research area, explaining the main ideas behind the recent progress and some of the  important challenges that lie ahead. 	
	\end{abstract}
	
	
	\thanks{The  author was supported  by the  Research Grant ELIDEK HFRI-NextGenerationEU-15689.}
	
	
	
	\date{\today}

	\maketitle
	
	\setcounter{tocdepth}{1}
	\tableofcontents

\section{Introduction}

A general principle in arithmetic Ramsey theory is that if you finitely partition
$\N:=\{1,2,\ldots\}$, then a cell of the partition will have a lot of structure. By structure, we mean solutions to algebraic equations or systems of algebraic equations.
This is nicely illustrated by a theorem of  Schur~\cite{Sc16} and van der Waerden~\cite{W27} and its extension by Rado~\cite{R33}, which covers systems of linear equations.
In the case of a single linear equation, it takes a particularly nice form that we will explain next.
\begin{definition}
	We say that the equation $P(x,y,z)=0$ is {\em partition regular}  if for every finite partition of $\N$ there exist distinct $x,y,z$ on the same cell that satisfy the equation.
\end{definition}
Rado's theorem states that for $a,b,c\in \N$, the equation
$$
ax+by=cz
$$
is partition regular if and only if  either $a,b,$ or $a+b$ equals $c$; we call
such a triple of positive integers $(a,b,c)$  a {\em Rado triple}.  To prove the  necessity of the conditions, take a prime $p>a+b+c$ and partition the positive integers into $p-1$ cells according to the last non-zero digit in their base $p$ expansion. Such colorings
arise from   level sets of some modified Dirichlet characters  (see Section~\ref{SS:pretex}) and will be our main source of counterexamples in this article.

Higher degree polynomial equations have proved much harder to tackle.
Perhaps the most well-known problem in this area is that of Erd\H os and Graham \cite{Gr07,Gr08}, which asks whether Pythagorean triples are partition regular, i.e. whether the equation
\begin{equation}\label{E:pyth}
x^2+y^2=z^2
\end{equation}
is partition regular.
Graham in \cite{Gr07} places the origin of the problem in the late 70's and offered \$250 for its solution, noting that ``there is actually very little data (in either direction) to know which way to guess.''  While this may have been true a decade ago,
there have been some positive developments in recent years.
The case where one allows only partitions of $\N$ into two sets was verified in 2016 with the help of a computer search \cite{HKM16}; this endeavor was hailed as the ``longest mathematical proof'' at the time, consuming 200 terabytes of data \cite{L16}. Also, partition regularity of Pythagorean triples  for level sets of finite-valued multiplicative functions  was recently established \cite{FKM23}.

The seminal work of  Furstenberg \cite{Fust} and S\'ark\"ozy \cite{Sa78-1}, culminating in the influential polynomial Szemer\'edi theorem of Bergelson and Leibman \cite{BL96}, allows us to deal with  shift-invariant systems of polynomial equations.  Regarding    equations that resemble \eqref{E:pyth},  here are some   examples  that are known to be or not to be partition regular:
\begin{enumerate}
	\item  $x^2+y=z$ is partition regular (Bergelson~\cite{Be87});
	
	\item $x+y=z^2$ is not partition regular (Csikv\'ari-Gyarmati-S\'ark\"ozy~\cite{CGS12});
	
	\item   $x^2+y^2=z$ is not partition regular (variant of Green-Lindqvist~\cite[Section~2]{GL19});
	
	\item   $x^2+y=z^2$ is  partition regular (Moreira~\cite{Mo17}).
\end{enumerate}
 The second equation is  partition regular if we consider partitions with only two cells~\cite{GL19,P18} or if we only require $x,y$ to belong to the same cell~\cite{KS06}.
There are a variety of tools used to prove these results, such as elementary combinatorial techniques, Fourier analysis, topological dynamics, and ergodic theory.
Other partition regularity results  of similar flavor can be found in \cite{Al22,BLM21, BJM17,DL18,DR18,DLMS23, L91}.

However, in the fully non-linear setting, the only known results involve equations in more than three variables. For example,  a
result of  Chow-Lindqvist-Prendiville \cite{CLP21} establishes that the equation
$$
x_1^2+x_2^2+x_3^2+x_4^2=x_5^2
$$
 is partition regular (see also \cite{BP17,Ch22, CC24,  P21} for  related results).

Here we focus on homogeneous quadratic equations in three variables, for which the methodology developed to deal with previous results has so far proved inadequate. A problem that drives much recent research in this area is the following:
\begin{conjecture}\label{Conj1}
If  $a,b,c\in \N$, then the equation
	\begin{equation}\label{E:pythgen}
		ax^2+by^2=cz^2
	\end{equation}
is partition regular if and only if $(a,b,c)$ is a Rado triple.
\end{conjecture}
The necessity of the conditions follows by combining  Rado's theorem \cite{R33} with the fact that  partition regularity of the equation  $P(x^2,y^2,z^2)=0$ implies partition regularity of the equation $P(x,y,z)=0$.
Regarding sufficiency, unfortunately there is no instance of a Rado triple for which the partition regularity of \eqref{E:pythgen} has been proved.
Three typical cases to keep in mind are the equations
 $$
 x^2+y^2=z^2,\quad  x^2+2y^2=z^2, \quad x^2+y^2=2z^2,
 $$
 which, as we will see, have particular difficulties that  are representative of the general case. The last equation was conjectured to be partition regular by  Gyarmati-Ruzsa~\cite{GR} and may even be density regular with respect to the natural density, i.e. it may have non-trivial solutions on any subset of the integers with positive upper density.

 Although we are not in a position to address the question of partition
regularity for triples, we will describe a systematic method, based on recent work by the author, Klurman, and Moreira~\cite{FKM23,FKM24}, and older work by the author and Host~\cite{FH17}, which allows us to address the
problem for pairs in several cases.
\begin{definition}
	Given non-zero $a,b,c\in \Z$, we say that the equation \eqref{E:pythgen} is {\em partition regular with respect to $x,y$}  if for every finite coloring of $\N$ there exist distinct $x,y\in \N$, with the same color, and $z\in \N$ that satisfy \eqref{E:pythgen}. Similarly, we define the partition regularity of \eqref{E:pythgen}  with respect to the variables $x,z$ and $y,z$.
\end{definition}
One of our main results in \cite{FKM23} is that Pythagorean pairs are partition regular, which answers a question that
 appeared in  \cite{DLMS23, Fr16, FH17, G24}:
\begin{theorem}
	The Pythagorean equation \eqref{E:pyth} is partition regular with respect to all pairs of variables.
	\end{theorem}
	Our starting point  is to
	link our problem with one concerning asymptotic properties of bounded  multiplicative functions
	via a representation result of Bochner-Herglotz;
	 the details of this link are explained in Section~\ref{S:proof}. Let us just say that multiplicative functions play the role of characters in multiplicative Fourier analysis, which is an alternative to  additive Fourier analysis that is worth exploring for the Pythagorean equation, or any other equation that is dilation-invariant but not translation-invariant.
As we explain in Section~\ref{SS:multi} below, there is no known analogous   ``Bochner-Herglotz representation result''  in the context of the partition regularity of Pythagorean triples, and an approach based on ergodic theory may be better suited to this problem. This is similar to the situation with shift-invariant patterns, where patterns of length two can be handled by a variety of tools, including Fourier analysis, but ergodic theory tools are better suited for patterns of length three and above, especially when these patterns are sparse.

In \cite{FKM24}, we extended our methodology to cover  partition regularity with respect to the variables $x,y$ of  more general equations of the form $ax^2+by^2=cz^2$
 where $a,b,c$ are non-zero integers such that either $ac$ or $bc$ is a square. We will discuss why equations like  $x^2\pm y^2=2z^2$ are harder to handle.
In all of these problems, our approach depends on an in-depth study of the asymptotic behavior of the averages
$$
\frac{1}{N^2}\sum_{1\leq m,n\leq N} f(P_1(m,n))\cdot \overline{f(P_2(m,n))},
$$
where $f$ is an arbitrary completely multiplicative function taking values on the unit circle, and $P_1,P_2$ are binary quadratic forms  used to parametrize the
solutions of \eqref{E:pythgen}.  We use Gowers-uniformity properties to deal with ``minor-arc'' multiplicative functions; this is the aperiodic case discussed in Section~\ref{S:aperiodic}, and concentration estimates to deal with ``major-arc'' multiplicative functions; this is the pretentious case discussed in Section~\ref{S:concentration}.

Before these works, partition regularity of pairs was proved by the author and Host in \cite{FH17}, which covered, for example,  partition regularity of the equation $16x^2+9y^2=z^2$ with respect to $x,y$, but missed the case of Pythagorean pairs for reasons we will explain later.
Extending these ideas, Sun in \cite{Su18,Su23} established partition regularity  with respect to  $x,y$ for the equation $x^2-y^2=z^2$ when $\N$ is replaced by the ring of integers of a larger number field, such as the Gaussian integers. However, the methods used there do not apply to $\N$.

One may also pursue similar partition regularity questions on the rationals, but as it turns out, in the case of homogeneous equations (or dilation invariant patterns), partition regularity over $\Q$ and $\Z$ are equivalent problems~\cite[Proposition~2.13]{DLMS23}. This is not the case
for non-homogeneous patterns, and there are several interesting patterns that are known to be partition regular over the rationals, but it is not known whether they are partition regular over the integers, see for example some recent results in \cite{A23,BS24}.
\subsection*{Notation} \label{SS:notation}
We let  $\N:=\{1,2,\ldots\}$,  $\Z_+:=\{0,1,2,\ldots \}$,
 $\R_+:=[0,+\infty)$, $\S^1\subset\C$ be the unit circle, and $\U$ be  the closed complex unit disk.
With $\P$ we
denote the set of primes and throughout we use the letter $p$ to denote prime numbers.

For $t\in \R$, we let $e(t):=e^{2\pi i t}$.
For $z\in \C$, with $\Re(z)$, $\Im(z)$  we denote the real  and imaginary parts of $z$ respectively, and we also let  $\exp(z):=e^z$.

For  $N\in\N$, we let $[N]:=\{1,\dots,N\}$. We often denote sequences $a\colon \N\to \U$
by  $(a(n))$, instead of $(a(n))_{n\in\N}$.

If $A$ is a finite non-empty subset of the integers and $a\colon A\to \C$, we let
$$
\E_{n\in A}\, a(n):=\frac{1}{|A|}\sum_{n\in A}\, a(n).
$$
We write $a(n)\ll b(n)$ if for some $C>0$ we have $a(n)\leq C\, b(n)$ for every $n\in \N$.

Throughout this article, the letter $f$ is typically  used for multiplicative functions  and the letter $\chi$ for Dirichlet characters. With $\CM$ we denote the set of completely multiplicative functions with modulus $1$.

\subsection*{Acknowledgement}
I  am grateful   to my  co-authors  on the articles  \cite{FH17,FKM23,FKM24},  B.~Host, O.~Klurman, and J.~Moreira;  without them,
 these works would not have been possible. I would also like to thank A.~Mountakis and D.~Sclosa for a number of  useful comments.

\section{Partition regularity  of generalized Pythagorean pairs}
We state here the main results concerning the partition regularity of pairs for homogeneous quadratic equations, and also make a crucial reduction to a density regularity statement involving parameterizations of the solution sets of the respective pairs.

\subsection{Main results}
A variant of Conjecture~\ref{Conj1} that covers partition regularity of pairs  for \eqref{E:pythgen} is the following one:
\begin{conjecture}\label{Conj2}
	If $a,b,c\in \Z$ are non-zero and
	at least one of the integers $ac,bc, (a+b)c$ is a square, then \eqref{E:pythgen}  is partition regular with respect to $x,y$.
\end{conjecture}
If \eqref{E:pythgen} is partition regular with respect to all three variables, then it is obviously partition regular for any pair of variables.
However, the converse is not true: \cite[Theorem 1.1]{FKM23} implies that the equation $x^2+y^2=4z^2$ is partition regular with respect to any pair of variables, but this equation is not partition regular with respect to all three variables since $(1,1,4)$ is not a Rado triple.

In  \cite{FKM24}  we prove the following result, which partially answers  Conjecture~\ref{Conj2}.
	\begin{theorem}\label{T:PairsPartition}	
		If  $a,b,c\in \Z$ are non-zero and  $ac$ or $bc$ is a square, then \eqref{E:pythgen} is partition regular with respect to $x,y$.
\end{theorem}
Taking $a=b=c=1$ and $a=-b=c=1$  we get that the equation $x^2+y^2=z^2$ is partition regular with respect to all pairs of variables $x,y,z$. Taking $a=c=1, b=2$, and $a=c=1, b=-2$ we get that the equation
$x^2+2y^2=z^2$ is partition regular with respect to the pairs  $x,y$ and $y,z$. We will see later
why the pair $x,z$ is harder to handle for this equation, and why a similar problem arises with the equation
$x^2+y^2=2z^2$, this time with respect to all pairs of variables.
In this regard, we can only give a conditional result, assuming a difficult
Chowla-Elliott type conjecture (see Section~\ref{SS:apChowla}) for a class of ``sufficiently aperiodic'' completely multiplicative functions. 
	\begin{theorem}\label{T:PairsPartitionCond}
	Suppose that Conjecture~\ref{C:2quadratic} below holds.
		If  $a,b,c\in \Z$  are non-zero and  $(a+b)c$ is a square (perhaps zero), then \eqref{E:pythgen} is partition regular with respect to $x,y$.
	\end{theorem}
Although the conditionality of this result is an unpleasant feature, we can cover, unconditionally,
 partitions generated by pretentious multiplicative functions (see \cite[Theorem~1.4]{FKM24} for the precise statement), which are the only known sources of counterexamples to the
 partition regularity of pairs of homogeneous quadratic equations.

An immediate consequence of Theorems~\ref{T:PairsPartition} and \ref{T:PairsPartitionCond} is that
conditionally on the truth of Conjecture~\ref{C:2quadratic} in Section~\ref{SS:apChowla} we have
that if $(a,b,c)$ is a Rado triple, then  the equation \eqref{E:pythgen} is partition regular with respect to all pairs of variables $x,y,z$.
We get this result unconditionally if $a=c$ or $b=c$.

\subsection{Parametric reformulation}\label{SS:parametric}
To prove our main results we will work with the  parametric form of the solutions of equation \eqref{E:pythgen}. We give three representative examples of such parametrizations
(below, $k,m,n$ are arbitrary integers):
\begin{enumerate}
	\item \label{PM1} $x^2+y^2=z^2$:
	$x=k\, (m^2-n^2)$,   $y=k\, (2mn)$,   $z=k\, (m^2+n^2)$;

\medskip
	
	\item \label{PM2} $x^2+2y^2=z^2$:	
	$x=k\, (m^2-2n^2)$,   $y=k\, (2mn)$,   $z=k\, (m^2+2n^2)$;
	
\medskip
	
	\item \label{PM3} $x^2+y^2=2z^2$:	
	$x=k\, (m^2-n^2+2mn)$,   $y=k\, (m^2-n^2-2mn)$,   $z=k\, (m^2+n^2)$.
\end{enumerate}

To give the reader a sense of how  certain features of these parametrizations affect our analysis, we make a few remarks
that  will be properly justified when we go deeper into the proof outline of our main results.

In the case of \eqref{PM1}, in order to prove partition regularity with respect to $x,y,$ it suffices to show that for every finite partition of $\N$ there exist $k,m,n\in \N$ such that the integers
$$
k(m^2-n^2) \quad \text{and} \quad k(2mn)
$$
are distinct and belong to the same partition cell.
The fact that both polynomials $m^2-n^2$ and $2mn$  are products of linear forms simplifies the argument substantially.   On the other hand, the parameter $z$ uses the irreducible quadratic form $m^2+n^2$, which takes values on a zero density subset of $\N$,  and this makes the corresponding analysis of partition regularity with respect to the variables $x,z$ (or $y,z$) much more complicated.  On the positive side, the polynomial $m^2+n^2$ is a norm form  (higher degree irreducible polynomials that do not have this property are out of reach at the moment) and the corresponding ring of integers one needs to work with in part of our analysis is the Gaussian integers, which is a principal ideal domain with finitely many units. All these things reduce the complexity of our argument.

In the case of \eqref{PM2}, handling the parameters $y,z$ is similar to the case of \eqref{PM1}. However, handling the parameters $x,y$ is much harder, since the ring of integers associated with the polynomial $m^2-2n^2$ is $\Z[\sqrt{2}]$, which has infinitely many units; this adds another layer of problems to overcome. Finally, the variables $x,z$ depend on two irreducible
quadratic forms, which is the most difficult case to handle, and we can offer only conditional results in this case.
The same problem appears in the equation \eqref{PM3}, in this case all pairs depend on two irreducible quadratic forms, so  we can only prove conditional results in these cases.

\subsection{Density regularity}
It will be more convenient for us
 to prove stronger versions of   Theorems~\ref{T:PairsPartition} and \ref{T:PairsPartitionCond}, which  deal with density regularity. Roughly, the stronger statements will state that every set with positive multiplicative density contains the pairs we are interested in, written in parametric form.

 We recall some standard notions.
 A \emph{multiplicative F\o lner sequence in $\N$} is a sequence $\Phi=(\Phi_K)_{K=1}^\infty$ of finite subsets of $\N$ that is  asymptotically invariant under dilation, in the sense that
 $$
 \lim_{K\to\infty}\frac{\big|\Phi_K \cap (x\cdot \Phi_K)\big|}{|\Phi_K|}=1 \quad \text{for all } x\in \N.
 $$
 An example of a multiplicative F\o lner sequence is given by
 $$
 \Phi_K	:=\Big\{\prod_{p\leq K}p^{a_p}\colon 0\leq  a_p\leq K\Big\}, \quad K\in \N.
 $$	
 The \emph{upper multiplicative density} of a set $\Lambda\subset\N$ with respect to a multiplicative F\o lner sequence $\Phi=(\Phi_K)_{K=1}^\infty$ is the quantity
 $$\bar{d}_\Phi(\Lambda):=\limsup_{K\to\infty}\frac{\big|\Phi_K
 	\cap \Lambda\big|}{|\Phi_K|}.$$
 		For example, for every $r\in \N$, the set $r\Z+j$ has multiplicative density $1$ if $j=0$ and multiplicative density $0$ if $j=1,\ldots, r-1$, so a set with positive multiplicative density should contain many elements that are highly divisible.    The set of squares has multiplicative density $0$, while  the set of integers $n\in \N$ such that $2^{2k}\mid n$ but $2^{2k+1}\nmid n$ for some $k\in \Z_+$,  has multiplicative density $1/2$.
\begin{definition}
	We say that the pair  $P_1,P_2\in \Z[m,n]$ is  {\em good for density regularity} if  whenever $\Lambda\subset \N$ satisfies  $\bar{d}_\Phi(\Lambda)>0$
	for some multiplicative F\o lner sequence $\Phi$,  there exist  $m,n\in\N$, such that
	$P_1(m,n)$ and  $P_2(m,n)$ are distinct positive integers and
	\begin{equation}\label{E:dpos}
	\bar{d}_\Phi\big( (P_1(m,n))^{-1}\Lambda\cap (P_2(m,n))^{-1}\Lambda\big)>0,
	\end{equation}
where $r^{-1}\Lambda:=\{n\in \N\colon rn\in \Lambda\}$.
\end{definition}
If the pair $P_1,P_2\in \Z[m,n]$ is good for density regularity, then for any $\Lambda\subset \N$ that satisfies  $\bar{d}_\Phi(\Lambda)>0$ for some multiplicative F\o lner sequence $\Phi$, the intersection $(P_1(m,n))^{-1}\Lambda\cap (P_2(m,n))^{-1}\Lambda$ has positive multiplicative density, and in particular is non-empty for some $m,n\in\N$.
Taking any $k\in(P_1(m,n))^{-1}\Lambda\cap (P_2(m,n))^{-1}\Lambda$, it follows that  $kP_1(m,n)$, $kP_2(m,n)$ are distinct and they both belong to $\Lambda$.

By considering appropriate parametrizations of solutions of \eqref{E:pythgen} when $ac$ or $bc$ is a square, it turns out that in order to prove Theorem~\ref{T:PairsPartition} it suffices to establish the following result:
\begin{theorem}\label{T:PairsDensityParametric}
	Suppose that for some $\alpha,\gamma\in \N$  and non-zero $\beta \in \Z$
	we have
	$$
	P_1(m,n)=\alpha m^2 +\beta n^2, \qquad P_2(m,n)=\gamma  mn.
	$$
	Then the pair of polynomials $P_1,P_2$ is good for density regularity.
\end{theorem}
More generally, our argument shows density regularity if at least  one of the two polynomials is reducible and not a square, and the two polynomials are not multiples of each other.

On the other hand,  appropriate parametrizations of solutions of \eqref{E:pythgen} when $(a+b)c$ is a square show that
Theorem~\ref{T:PairsPartitionCond} is an immediate consequence of  the following result:
\begin{theorem}\label{T:PairsDensityParametricIrreducible}
	Suppose that for some $\alpha,\beta,\alpha',\beta'\in \Z$
	the  polynomials
	\begin{equation}\label{E:2irre}
	P_1(m,n)=m^2+\alpha mn+\beta n^2, \qquad P_2(m,n)=m^2+\alpha' mn+\beta' n^2
	\end{equation}
	are distinct and irreducible, and  Conjecture~\ref{C:2quadratic} holds
	for $P_1, P_2$ (see Section~\ref{SS:apChowla}).
	Then the pair of polynomials $P_1,P_2$ is good for density regularity.
\end{theorem}
More generally, our argument shows density regularity, if
$P_1,P_2$ are irreducible quadratic forms, not multiples of each other, such that either $P_1(1,0)=P_2(1,0)$, or $P_1(0,1)=P_2(0,1)$.  If neither of the last two conditions on  the coefficients of $m^2$ and $n^2$  is met, then the  corresponding pair of polynomials is often not partition regular.

We note that our arguments also give that the positivity property \eqref{E:dpos} occurs quite often, in fact it holds for a set of $m,n\in \N$ with positive lower density
(taken over the squares $[N]\times [N]$).

\subsection{Other results} Our methodology is quite flexible and can be applied to other partition regularity problems of
pairs of dilation invariant patterns and related results. Here is a list of some additional examples:
\begin{enumerate}
\item (More general homogeneous quadratic equations)
The methods of this article can be used to prove partition and density regularity of more general equations $P(x,y,z)=0$, where $P$ is a homogeneous quadratic form.
For example, take  $a,b\in \N$, $d\in \Z$, and consider the equation
\begin{equation}\label{E:xyz}
ax^2+by^2+dxy=az^2.
\end{equation}
Direct computation shows that \eqref{E:xyz}  is satisfied when
$$
x=k(bm^2-an^2), \quad y=km(-dm+2an), \quad z=k(bm^2-dmn+an^2),
$$
where  $k,m,n\in\Z$.
Using the remark following Theorem~\ref{T:PairsDensityParametric}, we get that \eqref{E:xyz} is density regular with respect to  $x,y$ and $y,z$. We also get density regularity with respect to $x,z$ if
the polynomial $bm^2-an^2$ is reducible, i.e. if $ab$ is a square. The last assumption can be removed if we assume that  Conjecture~\ref{C:2quadratic} in Section~\ref{SS:apChowla} holds; to see this,   use the remark
 following Theorem~\ref{T:PairsDensityParametricIrreducible}.

\smallskip

 \item (Pythagorean triples on level sets of multiplicative functions) In \cite[Theorem~1.7]{FKM23} we showed that for any completely multiplicative function $f\colon \N\to \S^1$ taking finitely many values, there exist $x,y,z\in \N$ such that
  $$
  x^2+y^2=z^2 \quad \text{and} \quad  f(x)=f(y)=f(z)=1.
  $$

\smallskip

\item (A question from \cite{DLMS23})
Answering a question of Donoso-Le-Moreira-Sun from \cite{DLMS23} we showed in \cite[Theorem~1.8]{FKM23} (our argument also uses the main result from \cite{Tao15})  that the patterns $(km^2, k(n^2+n))$ are density regular.

\smallskip

\item (An extension of a result  of Klurman-Mangerel)
Recently, Charamaras-Mountakis-Tsinas~\cite{CMT24} used a variant of the
$Q$-trick, which we explain in Section~\ref{S:Qtrick}, to show, among other things, that for every  $a\in \N$  and every pair of completely multiplicative functions $f,g\colon \N\to \S^1$ we have
\begin{equation}\label{E:fabc}
\liminf_{n\to\infty}|f(an+1)-g(an)|=0,
\end{equation}
and we cannot replace  $1$ by $2$ even when $a=1$.
When $f=g$, equation~\eqref{E:fabc}  was previously established by Klurman-Mangerel~\cite[Theorem~1.1]{KM18} for  $a=1$ and for general $a\in \N$ by  Donoso-Le-Moreira-Sun~\cite[Corollary~1.7]{DLMS23}
\end{enumerate}

\section{A model problem and connection with multiplicative functions}

\subsection{An additive model-problem} Our approach to solving the density regularity problems of Theorems~\ref{T:PairsDensityParametric} and \ref{T:PairsDensityParametricIrreducible}
	 tries to model
	 Furstenberg's ergodic-theoretic proof of S\'ark\"ozy's theorem. Let us briefly describe this argument.  The goal is to show that if $\Lambda$ is a positive density subset of the integers, i.e.
	 $$
	 \overline{d}(\Lambda):=\limsup_{N\to\infty}\frac{|\Lambda\cap [N]|}{N}>0,
	 $$
	 then  there exists $n\in \N$ such that
	 \begin{equation}\label{E:denspos}
	 \overline{d}(\Lambda\cap(\Lambda-n^2))>0.
	 \end{equation}
	 Using a correspondence principle of Furstenberg~\cite{Fust}, it suffices to show that for every ``additive'' measure-preserving system $(X,\mu,T)$, i.e. probability space $(X,\mathcal{X},\mu)$, measure-preserving transformation $T\colon X\to X$,  and measurable set $A$ with $\mu(A)>0$,
	 there exists $n\in \N$ such that $\mu(A\cap T^{-n^2}A)>0$. In fact, we  show the stronger property
	 \begin{equation}\label{E:Sarkozy}
	 \liminf_{N\to\infty}\E_{n\in [N]}\, \mu(A\cap T^{-n^2}A)>0.
	 \end{equation}
	
	To prove \eqref{E:Sarkozy}, we have two ways to proceed. One is to first show that the rational Kronecker factor is characteristic (by appealing to the Hilbert space version of van der Corput's lemma), and then to deal with the rational Kronecker factor of the system using the good divisibility properties of $n^2$. Unfortunately, this approach cannot be modelled very well in our multiplicative setting, since our patterns use a mixture of addition and multiplication which does not allow cancellation when using a multiplicative version of van der Corput's lemma.

The second way to prove \eqref{E:Sarkozy} is to use a representation theorem for sequences of the form $n\mapsto \mu(A\cap T^{-n}A)$, or more generally
positive density sequences in $\Z$.
	 It follows from Heglotz's theorem that there exists a (positive) bounded Borel measure $\sigma$ on $[0,1)$  (depending on $A$) such that
	 \begin{equation}\label{E:spectraladditive}
	 \mu(A\cap T^{-n}A)=\int e(nt) \, d\sigma(t)
	 \end{equation}
	for every $n\in \Z$.  Averaging over $n\in [N]$,  taking the limit as $N\to \infty$, and using the ergodic theorem to lower bound the left side   and the bounded convergence theorem to compute the right side,  we deduce that $\sigma(\{0\})\geq (\mu(A))^2>0.$

	  It follows that in order to establish \eqref{E:Sarkozy} it suffices to show that if $\sigma$ is  a bounded (positive) measure on $[0,1)$ such that $\sigma(\{0\})=\delta>0$ and
	  $\widehat{\sigma}(n)\geq 0$ for every $n\in \N$, then
	  \begin{equation}\label{E:liminf1}
	  \liminf_{N\to\infty}\E_{n\in [N]}\, \int e(n^2t)\, d\sigma(t)>0.
	  \end{equation}
	  The advantage of this new setting is that exponential sums along squares
	  can be estimated efficiently, and this helps us analyze the left side in \eqref{E:liminf1}. Indeed, it is a well-known consequence of Weyl's theorem that if $t$ is irrational, we have
 \begin{equation}\label{E:additive0}
	   \lim_{N\to\infty}\E_{n\in [N]}\,  e( n^2t)=0.
	  \end{equation}
	   Using the bounded convergence theorem, this essentially allows us to reduce matters to handling the ``rational'' part of the measure $\sigma$, which can be done using the divisibility properties of $n^2$ as follows: For $A\subset [0,1),$ let
	  $$
	  \sigma_{A}:=\sigma|_{A}
	  $$
	 and  choose $Q\in \N$ so that
	  \begin{equation}\label{E:a1}
	  \sigma_{\Q}( [0,1)\setminus A_Q)\leq \delta^2/2  \quad \text{ where }  \quad
A_Q:=	  \{0,1/Q,\ldots, (Q-1)/Q\}.
	  \end{equation}
	  Using the positivity property  $\widehat{\sigma}(n)\geq 0, n\in \N$, we get that \eqref{E:liminf1} would  follow  from
	  \begin{equation}\label{E:a1.5}
	  \liminf_{N\to\infty}\E_{n\in [N]}\, \int e( Q^2n^2t)\, d\sigma(t)>0;
	  \end{equation}
	  indeed, the limit in $\eqref{E:liminf1}$ is at least $1/Q$ times the limit in \eqref{E:a1.5}.
	  Using \eqref{E:additive0} and the bounded convergence theorem,  we get that it suffices to show that
	      \begin{equation}\label{E:a2}
	  \liminf_{N\to\infty}\E_{n\in [N]}\, \int e(Q^2n^2t)\, d\sigma_{\Q}(t)>0.
	  \end{equation}
	Since $Q^2t^2=0 \pmod{1}$ for every $t\in A_Q$  we have
	        \begin{equation}\label{E:a3}
	  \int e( Q^2n^2t)\, d\sigma_{A_Q}(t)=\int 1\, d\sigma_{A_Q}(t)=\sigma(A_Q)\geq \sigma(\{0\})=\delta^2 \ \text{ for every } \   n\in \N.
		  \end{equation}
	  Combining  \eqref{E:a1} and \eqref{E:a3},  it follows that the limit in \eqref{E:a2} is at least $\delta^2/2$, and the proof is complete.

	  We also note that the use of ergodic theory in this argument is rather superficial,
	 since it was only employed to reduce the density positivity property  \eqref{E:denspos} to the analytic positivity property  \eqref{E:liminf1}. This step  could also be done directly by showing that
	  	the sequence $n\mapsto d^*(\Lambda\cap (\Lambda-n))$ is positive definite,
	  	where $d^*$ is a density  that agrees with $\overline{d}$ on $\Lambda$ and is defined along a subsequence of the integers so that all relevant limits exist.

	\subsection{Bird's eye view of the modified strategy}
		We plan to use the general principles of the strategy explained in the previous subsection
	to establish the density statements of Theorems~\ref{T:PairsDensityParametric} and \ref{T:PairsDensityParametricIrreducible}, with some necessary modifications that we will briefly summarize now. In the following sections, we will give a much more detailed account on how to implement these modifications
	and it will become clear why the exact form of the homogeneous polynomials we have to work with   greatly affects the difficulty of the problem.

	\begin{enumerate}
\item (Representation result)	We can replace the ``additive'' measure-preserving system $(X,\mu,T$) by a measure-preserving system $(X,\mu,T_n)$ that has  multiplicative structure, but we opt for a more  direct approach and represent directly some  multiplicative  densities that are naturally linked to our problem. The fact that these densities are dilation invariant  and not translation invariant  leads to a representation result analogous to \eqref{E:spectraladditive}, where  the place of the linear exponential sequences $n\mapsto e(nt)$ is taken by
sequences with multiplicative structure, i.e. sequences $f\colon \N\to \S^1$ that satisfy $f(mn)=f(m)\cdot f(n)$ for every $m,n\in \N$.  See \eqref{E:dPhi}  for the exact identity.

  \medskip

\item (Reduced analytic statement) Using the previous representation result,  we reduce the proof of Theorems~\ref{T:PairsDensityParametric} and \ref{T:PairsDensityParametricIrreducible}   to showing  a positivity property analogous to \eqref{E:liminf1}, which involves averages of integrals of multiplicative functions evaluated over homogeneous polynomials in two variables.  See Theorems~\ref{T:MainMulti1} and \ref{T:MainMulti2} below for the exact statements.
It is the case, however, that we have a much larger variety of multiplicative functions than we have linear exponential sequences, and studying their asymptotic behavior turns out to be a much more delicate and interesting task.
Similar to the additive case, we consider two complementary cases, the ``structured'' and  the  ``random'' cases,  to handle the ``major-arcs'' and ``minor-arcs''   of the spectral measure.

\medskip

\item 	(Minor arcs) In the additive setting, we used Weyl's theorem to show that for irrational numbers    we have the vanishing property
 \eqref{E:additive0}. In the multiplicative setting, the place of irrational numbers is  taken by    the ``aperiodic'' multiplicative functions, for which we prove, using a rather involved argument, several vanishing statements about their correlations, see Section~\ref{S:aperiodic} for details.

	\medskip

\item(Major arcs)  In the additive setting, we were able to handle the part of the spectral measure $\sigma$ that was supported on rational numbers    by passing to a subprogression and then  using  the identity \eqref{E:a3}. In the multiplicative setting, the place of rational numbers  is taken by the ``pretentious'' multiplicative functions for which we are able to prove concentration estimates as a substitute for the identity \eqref{E:a3}. The exact statements appear in Section~\ref{S:concentration}.

\medskip

\item (The $Q$-trick) In the additive setting, it was easy to combine the information about the major and minor arcs  to simplify the expression \eqref{E:a1.5} to a point where  the positivity property  is obvious.  In the multiplicative setting, the simpler expression we arrive at has no obvious positivity property, so  it will be much less clear how to conclude. A key manipulation here is  what we call  the $Q$-trick, which allows us to obtain a property analogous to \eqref{E:a2} for some highly divisible $Q$.
We explain in Section~\ref{S:Qtrick} how this key  maneuver is used to complete the proof. 	
	 \end{enumerate}

	\subsection{From partition regularity to multiplicative functions}\label{S:proof} The starting point   for  the proofs  of Theorems~\ref{T:PairsDensityParametric} and \ref{T:PairsDensityParametricIrreducible}
	 is to reformulate  them  as  positivity properties involving integrals of  averages of completely multiplicative functions. As we also mentioned before, this key step can only be performed for partition and density regularity results of pairs, and has no useful analog we know of for triples (though see Section~\ref{SS:multi} for an  ergodic reformulation).

\subsubsection{Multiplicative functions}	A function $f\colon \N\to \U$, where $\U$ is the complex unit disk,  is called {\em multiplicative} if
$$
f(mn)=f(m)\cdot f(n)  \quad \text {  whenever  }  (m,n)=1.
$$
If necessary, we extend $f$ to an even function in $\Z$ by letting $f(0):=0$ and $f(-n):=f(n)$ for $n\in \N$.
It is called {\em completely multiplicative} if the previous identity holds for all $m,n\in\N$.
We let
$$
\CM:=\{f\colon \N\to \S^1 \text{ is a completely multiplicative function}\},
$$
where $\S^1$ is the  unit circle.
Throughout, we assume that $\CM$ is equipped with the topology of pointwise convergence; then $\CM$ is a metrizable compact space with this topology.

	\subsubsection{Reduction to a statement about multiplicative functions}\label{SS:analyticstatement}
	We will use  a representation result  of Bochner-Herglotz for positive definite functions  on locally compact Abelian groups. We  need it only for the discrete group $\Q_+$ with multiplication, in which case
	we can identify $\CM$ with the  Pontryagin dual of $(\Q_+,\times)$.
		We say that  $A\colon \Q_+\to \C$  is positive definite on  $(\Q_+,\times)$ if
	 for every $n\in\N$, all $r_1,\dots,r_n\in\Q^+$ and all $\lambda_1,\dots,\lambda_n\in\C,$ we have
	$$
	\sum_{i,j=1}^n\lambda_i\overline{\lambda_j}\,A(r_i\,r_j\inv)\geq 0.
	$$
		\begin{theorem}[Bochner-Herglotz]\label{T:BH}
		Let $A\colon \Q_+\to \C$ be  positive definite  on $(\Q_+,\times)$. Then there exists a finite Borel measure $\sigma$ on $\CM$ such that
		$$
		A(r/s)=\int_\CM f(r)\cdot \overline{f(s)}\, d\sigma(f), \quad \text{for all } r,s\in \N.
		$$
	\end{theorem}
		To use this representation result in our setting, we can argue as follows (an alternative way is via ergodic theory).   Given a  F\o lner  sequence $\Phi=(\Phi_N)$  of subsets of $\N$
 we can pass to a subsequence $\Phi'=(\Phi'_N)$ along which  $d_{\Phi'}(\Lambda)=d_{\Phi}(\Lambda)$ and  the limits
	$$
	d_{\Phi'}\big((r^{-1}\Lambda)\cap(s^{-1}\Lambda)\big) :=
	\lim_{N\to\infty}\frac{\big|(r^{-1}\Lambda)\cap(s^{-1}\Lambda)\cap \Phi'_N|}{|\Phi'_N|}
		$$
exist for all $r,s\in \Z$. Then the  sequence $r/s \mapsto 	d_{\Phi'}\big((r^{-1}\Lambda)\cap(s^{-1}\Lambda)\big) $ from $(\Q_+,\times)$ to $ [0,1]$  is well defined (we use that $d_{\Phi'}$ is invariant under dilations), and  a direct computation shows that it is 	positive definite. Using Theorem~\ref{T:BH}, we get that there exists a finite Borel measure $\sigma$  on $\mathcal{M}$  such that
\begin{equation}\label{E:dPhi}
 	d_{\Phi'}\big((r^{-1}\Lambda)\cap(s^{-1}\Lambda)\big)=
 	\int_\CM f(r)\cdot \overline{f(s)}\, d\sigma(f), \quad \text{for all } r,s\in \N.
\end{equation}
Furthermore, it is not hard to show (see \cite[Section~10.2]{FH17}) that
$$
\sigma(\{1\})\geq (d_{\Phi'}(\Lambda))^2.
$$
Using the previous two properties, it is easy to see  that in order to prove Theorem~\ref{T:PairsDensityParametric},   it suffices to prove the following result:
	\begin{theorem}\label{T:MainMulti1}
		Let $\sigma$ be a  bounded Borel measure on $\CM$ such that $\sigma(\{1\})>0$ and
		\begin{equation}\label{E:positivers}
			\int_\CM f(r)\cdot \overline{f(s)}\, d\sigma(f)\geq 0\quad  \text{for every } r,s\in \N.
		\end{equation}
Let also 	$\alpha,\gamma\in \N$  and non-zero $\beta \in \Z$.	Then
		\begin{equation}\label{E:sigmapositive1}
			\liminf_{N\to\infty} \E_{m,n\in[N]}\,   \int_{\CM} f\big(\alpha m^2 +\beta n^2\big)\cdot \overline{f\big(\gamma  mn\big)}\, d\sigma(f)>0.
		\end{equation}
	\end{theorem}
In order to prove	
	Theorem~\ref{T:PairsDensityParametricIrreducible}, it suffices to establish the following result:
 	\begin{theorem}\label{T:MainMulti2}
 	Let $\sigma$ be as in Theorem~\ref{T:MainMulti1},
and   $\alpha,\beta,\alpha',\beta'\in \Z$
 be such that
 \begin{equation}\label{E:2irre}
 	P_1(m,n)=m^2+\alpha mn+\beta n^2, \qquad P_2(m,n)=m^2+\alpha' mn+\beta' n^2
 \end{equation}
 are distinct and irreducible quadratic forms. Suppose that either
 Conjecture~\ref{C:2quadratic}  in Section~\ref{SS:apChowla} holds
 for $P_1, P_2$, or the measure $\sigma$ is supported on the class of pretentious multiplicative functions (see Section~\ref{SS:pretentious}).	Then
 	\begin{equation}\label{E:sigmapositive2}
 		\liminf_{N\to\infty} \E_{m,n\in[N]}\,   \int_{\CM} f\big(P_1(m,n)\big)\cdot \overline{f\big(P_2(m,n)\big)}\, d\sigma(f)>0.
 	\end{equation}	
 \end{theorem}
	 Both results were established in \cite{FKM24}, with special cases previously obtained in \cite{FH17, FKM23}.
	
 One may  wonder about our choice of averaging over $m,n\in \N$ in \eqref{E:sigmapositive1} and \eqref{E:sigmapositive2}. Why not take iterated limits $\liminf_{N\to\infty}\E_{n\in[N]}\liminf_{M\to\infty}\E_{m\in[M]}$, or perhaps take  multiplicative averages over $m,n$?  Taking iterated limits would lead to much harder problems for ``minor arc'' multiplicative functions. Although in some cases the required vanishing property can be obtained for logarithmic averages using the work of Tao-Ter\"av\"ainen~\cite{TT19}, the exceptional class of ``minor-arc'' multiplicative functions for which this method is inapplicable would render this argument as not very useful. Finally,
 taking multiplicative averages over the parameters $m,n$ leads to problems that we are not able to handle,  both for ``major-arc'' and ``minor-arc'' multiplicative functions.

\subsection{Proof strategy} In order to establish the necessary positivity properties of Theorems~\ref{T:MainMulti1} and \ref{T:MainMulti2}, we  use  the theory of completely multiplicative functions and proceed in two independent steps, each using a distinctly different methodology.

\medskip

\begin{enumerate}
	\item (Minor arcs)
We first study the part of
 $\sigma$  that is supported on aperiodic multiplicative functions.  The key step here is to establish
  Gowers uniformity properties of arbitrary aperiodic multiplicative functions, and then using  them show  that the averages  in \eqref{E:sigmapositive1} vanish (this last step is trickier when $P(m,n):=\alpha m^2+\beta n^2$ is irreducible).
To prove the uniformity properties we combine the inverse theorem of Green-Tao-Ziegler~\cite{GTZ12}
with some   quite delicate equidistribution results of variable nilsequences.
 Our approach is described in more detail in  Section~\ref{S:aperiodic}.

\medskip

 \item (Major arcs)  Next, we study the complementary case, covering the part of $\sigma$ that is supported on pretentious multiplicative functions (see Definition~\ref{D:Pretentious}). In this case we can establish
 concentration estimates, a feature of pretentious multiplicative functions that is not shared by the aperiodic ones. The details are given in Section~\ref{S:concentration}.  We use these concentration estimates to  simplify substantially the averages in \eqref{E:sigmapositive1}  when the range of $m,n$ is restricted to a suitable positive density subset of  $\N\times \N$  depending on some parameter $Q$. We   then show that a suitable multiplicative average, over the parameter $Q$, of these simplified expressions is zero for all but an exceptional set of multiplicative functions and it is positive for the rest. We call this the $Q$-trick, and it easily implies the required positivity.  The details are given in Section~\ref{S:Qtrick}.
\end{enumerate}

	\section{Aperiodic multiplicative functions - Vanishing of correlations } \label{S:aperiodic}

\subsection{Definition and characterization} The first class of multiplicative functions that will play a crucial role in our arguments are the aperiodic ones.
\begin{definition}
	The multiplicative function $f\colon \N\to \U$ is {\em aperiodic} if for every $a,b\in \N$ we have
	\begin{equation}\label{E:aperiodic}
	\lim_{N\to\infty} \E_{n\in[N]}\, f(an+b)=0.
	\end{equation}
\end{definition}
Sometimes in the literature, the term {\em non-pretentious}
is used instead of the term aperiodic.
It can be shown that $f\colon \N\to \U$ is aperiodic if and only if
	\begin{equation}\label{E:AperiodicChi}
\lim_{N\to\infty}\E_{n\in[N]} \, f(n)\cdot \chi(n)=0 \quad \text{for every Dirichlet character } \chi .
\end{equation}

Examples of multiplicative functions that are aperiodic are the Liouville and the M\"obius function, this follows from a variant of the prime number theorem on arithmetic progressions.

A rather amazing and very useful fact for our purposes is that correlations of  aperiodic multiplicative functions have very strong vanishing properties, which we will describe in the following subsections.

\subsection{Gowers uniformity and inverse theorem}
We will define some seminorms on $\ell^\infty(\N)$ that will be crucial for our study. Roughly, they are given by  the limit as $N\to\infty$ of the
Gowers norms of the sequence, restricted to the interval $[N]$.

\begin{definition}[Gowers norms on $\Z_N$~\cite{G01}]
	Let $N\in \N$  and $a\colon \Z_N\to \C$. For $s\in \N$ the \emph{Gowers $U^s(\Z_N)$-norm} $\norm a_{U^s(\Z_N)}$  is defined inductively as follows:   For every $h\in\Z_N$ we write $a_h(n):=a(n+h)$. We let
	$$
	\norm a_{U^1(\Z_N)}:=|\E_{n\in\Z_N}a(n)|
	$$
	and for every $s\geq 1$ we let
	\begin{equation}
		\label{eq:def-gowers}
		\norm a_{U^{s+1}(\Z_N)}:=\Bigl(\E_{h\in\Z_N}\norm{a\cdot \overline a_h}_{U^s(\Z_N)}^{2^s}\Bigr)^{1/2^{s+1}}.
	\end{equation}
\end{definition}
For example,
$$
\norm a_{U^2(\Z_N)}^4=\E_{n,h_1,h_2\in\Z_N}\,  a(n)\, \overline{a(n+h_1)}\, \overline{a(n+h_2)}\, a(n+h_1+h_2)
$$
and for the $U^s(\Z_N)$ norms for $s\geq 3$ a similar but more complicated closed formula can be given.

\begin{definition}[Gowers semi-norms on $\N$]
If $a\colon \N\to \C$ is a bounded sequence, we let
  $$
  \norm{a}_{U^s(\N)}=\limsup_{N\to \infty}\norm{a_N}_{U^s(\Z_N)},
  $$
   where $a_N\colon \Z_N\to \C$ is the periodic extension of $a\cdot {\bf 1}_{[N]}$ to $\Z_N$.
	\end{definition}
 It can be shown that $\norm{\cdot}_{U^s(\N)}$ satisfies the triangle inequality, thus defining a seminorm on $\ell^\infty(\N)$, and for every $s\in \N$ we have
\begin{equation}\label{E:UkIncreases}
	\norm a_{U^s(\N)}\leq\norm a_{U^{s+1}(\N)}.
\end{equation}
We say that the bounded sequence $a\colon \N\to \C$ is $U^s$-uniform
if
$$
\norm{a}_{U^s(\N)}=0.
$$
For example, if $\alpha$ is an irrational number and $s\in \N$, it can be shown that  the sequence $(e(n^s\alpha))$ is  $U^s$-uniform but not $U^{s+1}$-uniform, while if $\alpha$ is a rational non-integer it is $U^1$-uniform but not $U^2$-uniform.
It can also be shown that if $c$ is a non-integer positive real number, then the
sequence $(e(n^c))$ is $U^s$-uniform for all $s\in \N$.

Using Fourier inversion and the Parseval identity in $\Z_N$, for $N\in \N$,  it can be shown that if
$a\colon \N\to \U$ is a sequence, then
$$
	\norm a_{U^2(\N)}
	\leq \limsup_{N\to\infty} \max_{\xi\in\Z_N}|\wh{a_N}(\xi)|^{\frac{1}{2}},
$$
where the Fourier transform is taken in $\Z_N$.
It is easy to deduce from this the equivalence
\begin{equation}
	\label{E:inverse2}
	\norm a_{U^2(\N)}=0 \quad \Longleftrightarrow \quad\lim_{N\to\infty} \sup_{\alpha\in \R} \big|\E_{n\in[N]}\,  a(n)\cdot e(n\alpha)\big|=0.
\end{equation}
In order to give a generalization of this result that covers $U^s$-uniformity for $s\geq 3$,
we have to first account for additional obstructions to uniformity. For example,  when $s=3$ it can be seen that
\begin{equation}\label{E:U3examples}
\norm a_{U^3(\N)}>0 \quad \text{when} \quad  a(n)=e(n^2\alpha+n\beta) \quad \text{or} \quad
a(n)=e( [n\alpha]n\beta),
\end{equation}
 where  $\alpha, \beta \in \R$.  These turn out to be   examples of  nilsequences, which  are defined as follows: If $X$ is an $s$-step nilmanifold, i.e.
 $X=G/\Gamma$, where $G$ is an $s$-step nilpotent Lie group and $\Gamma$ is a discrete cocompact subgroup, we let
 $$
 a(n)=F(g^n\cdot e_X), \quad n\in \N,
 $$
 where $g\in G$, $e_X:=\Gamma$, and $F\in C(X)$. Moreover, to cover the second example in \eqref{E:U3examples}, we must  allow a function $F$ that is  not continuous, but is  Riemann integrable in the orbit closure of $g$. It is an amazing fact that this class of sequences suffices to give  a generalization of \eqref{E:inverse2} that covers $U^s$-uniformity.
This is a consequence of the inverse theorem of Green-Tao-Ziegler~\cite{GTZ12}; first proved in \cite{GT08a} for $s=3$, which is the only case we will need for applications to partition regularity. We state it  in the following convenient form:
	\begin{theorem}\label{T:GTZ}
	Let $a\colon \N\to \U$ be a sequence and $s\in \N$. Then the following two properties are equivalent:
\begin{enumerate}
	\item
$
	\norm{a}_{U^{s+1}(\N)}=0;
	$

\item\label{I:inverse} For every $s$-step nilmanifold $X=G/\Gamma$ and  $F\in C(X)$ we have
	$$
	\lim_{N\to\infty}\sup_{g\in G}\big|\E_{n\in[N]}\,  a(n)\cdot F(g^n\cdot e_X)\big|=0.
	$$
\end{enumerate}	
\end{theorem}
Moreover, in checking the second condition, one can assume that $F$ is smooth, and
using what is known as the ``vertical Fourier decomposition'', one can further  assume
that $F$ is a non-trivial vertical nilcharacter, i.e., assuming that $G$ is $s$-step nilpotent but not $(s-1)$-step nilpotent, we have that there exists
$\chi\colon X_s\to \S^1$,  a non-trivial character
of the compact Abelian group $X_s:=G_s/(G_s\cap \Gamma)$, such that
\begin{equation}\label{E:vernil}
F(ux)=\chi(u)\cdot F(x)
\end{equation}
for every $x\in X$ and $u\in G_s$.  It easily follows from this that $\int F\, dm_X=0$.
We note  that in the case where $X=\T$,  non-trivial vertical nilcharacters are functions of the form  $F(x)=e( kx)$ for some non-zero $k\in \Z$.

We will see in the next subsection how one can use this result to establish the  $U^s$-uniformity of aperiodic multiplicative functions.

\subsection{Gowers uniformity of aperiodic multiplicative functions}
A very strong and perhaps a priori unexpected property, is that aperiodic multiplicative functions are Gowers uniform of all orders.

To verify that they are $U^2$-uniform
we will use the following  orthogonality criterion due to  Daboussi-K\'atai~\cite{DD82,K86}:
\begin{theorem}\label{T:DK}
	Let  $a\colon \N\to\U$  be a sequence such that
	\begin{equation}\label{E:DK}
		\lim_{N\to\infty} \E_{n\in[N]}\, a(pn)\cdot \overline{a(qn)}=0
	\end{equation}
	for all distinct primes $p,q$. Then  for every multiplicative function $f\colon \N\to \U$ we have
	$$
	\lim_{N\to\infty}\E_{n\in[N]}\,  f(n)\cdot a(n)=0.
	$$
\end{theorem}
Applying Theorem~\ref{T:DK} for   $a(n):= e( n\alpha)$,
 it follows  that
\begin{equation}\label{E:DKf}
\lim_{N\to \infty}\E_{n\in[N]}\,  f(n)\cdot e(n\alpha)=0
\end{equation}
for every irrational $\alpha$. Since aperiodic multiplicative functions, by definition, do not correlate with $ e(n\alpha)$ for $\alpha$ rational,
 it follows  that for those multiplicative functions  \eqref{E:DKf} holds for every $\alpha\in \R$. It is only slightly more difficult to show that a similar property holds even if we allow the phase $\alpha$ to depend on $N$~\cite[Lemma~9.1]{FH17}, namely
 \begin{equation}\label{E:DKuniform}
 	\lim_{N\to \infty}\sup_{\alpha\in\R}\big|\E_{n\in[N]}\,  f(n)\cdot e( n\alpha)\big|=0.
 \end{equation}
  A direct consequence of this identity and the characterisation of $U^2$-uniformity (see \eqref{E:inverse2}) is that
$$
\norm{f}_{U^2(\N)}=0.
$$
It is also true, but much harder to prove, that a similar property holds for higher order uniformity. This was first done for the Liouville and the M\"obius functions for $s=3$ by Green-Tao~\cite{GT08} and then
by Green-Tao-Ziegler~\cite{GT12,GTZ12} for general $s\in \N$. Their approach used a Type I  and Type II sums analysis to verify property~\eqref{I:inverse} of Theorem~\ref{T:GTZ}, and this required explicit quantitative control for averages of multiplicative functions along arithmetic progressions, not shared by  general aperiodic multiplicative functions (and the exceptional class is too large to be ignored in our applications). Using a somewhat different argument, the previous results were extended
  in a joint paper by the author and Host in \cite[Theorem~2.5]{FH17}  to general aperiodic multiplicative functions
  (see also  Matthiesen~\cite{M18} for an alternative approach of a related result).
\begin{theorem}\label{T:FH}
	Let $f\colon \N\to \U$ be an aperiodic  multiplicative function. Then for every $s\in \N$ we have
$$
\norm{f}_{U^s(\N)}=0.
$$
\end{theorem}
Let us  briefly sketch some key ideas of the proof of this result.
We start by using  Theorem~\ref{T:GTZ}. It reduces the problem to establishing  a strengthening of \eqref{E:DKuniform} where the place of $e(n\alpha)$ is taken by   nilsequences.
Namely, we need to show that if $X=G/\Gamma$ is a fixed $s$-step nilmanifold and $F\in C(X)$ is a  non-trivial vertical nilcharacter, then
\begin{equation}\label{E:nilkey}
\lim_{N\to \infty}\sup_{b\in G}\big|\E_{n\in[N]}\,  f(n)\cdot F(b^n\cdot e_X)\big|=0.
\end{equation}
Verifying this property is a rather laborious task, so let us just explain how it works in the much simpler setting where we do not take the sup over $b$, the nilmanifold $X$ is connected,  and the sequence  $(b^n\cdot e_X)$ is equidistributed on $X$.
So in this case, we use the above-mentioned orthogonality criterion of Daboussi-K\'atai
for $a(n):= F(b^n\cdot e_X)$. We get that it suffices to verify the vanishing property
\eqref{E:DK} for all distinct primes $p,q$, or equivalently, that
\begin{equation}\label{E:pq0}
\lim_{N\to\infty}\E_{n\in[N]}\, F(b^{pn}\cdot e_X)\cdot \overline{F(b^{qn}\cdot e_X)}=0.
\end{equation}
It can be shown that the sequence
$$
\big((b^{pn}\cdot e_X, b^{qn}\cdot e_X)\big)_{n\in\N}
$$
is equidistributed in a connected subnilmanifold $Y=H/\Delta$ of $X\times X$.
Equation~\eqref{E:pq0} would thus follow if we showed that
\begin{equation}\label{E:Y0}
\int F(x)\cdot  \overline{F(x')}\, dm_Y(x,x')=0.
\end{equation}
This is not a simple task, since the structure of the possible nilmanifolds $Y$ depends on $b$ in a way that is not easy to determine explicitly. The key property that allows us to prove \eqref{E:Y0} is that  no matter how complicated $Y$ may be it always satisfies the following
$$
\{(g^p,g^q)\colon g\in G\}\subset H\cdot (G_2\times G_2).
$$
This follows by using that the sequences $(b^{pn}\cdot e_X)$ and $(b^{qn}\cdot e_X)$ are both equidistributed in $X$. By taking iterated commutators of elements on the left we deduce the following key property
$$
(u^{p^s},u^{q^s})\in H \, \text{ for every } u \in G_s.
$$
This invariance property easily implies that $(x,x')\mapsto F(x)\cdot \overline{F(x')}$ is a vertical nilcharacter of $Y$ with non-zero frequency, so \eqref{E:Y0} holds.

If the place of the sequence $b^n$ takes a more general polynomial sequence, say $a^nb^{n^2}$ for some $a,b \in G$, we first establish that there exist normal subgroups $G^1,G^2$ of $G$  such that $G^1\cdot G^2=G$ and
$$
\{(g_1^pg_2^{p^2},g_1^qg_2^{q^2})\colon g_1\in G^1, g_2\in G^2\}\subset H\cdot (G_2\times G_2),
$$
and  by taking iterated commutators of elements on the left, we deduce   that
\begin{equation}\label{E:invariance}
\text{ the set} \quad   \{u\in G_s\colon (u^{p^j},u^{q^j})\in H \text{ for some } j\in \N\}\quad  \text{ generates } G_s.
\end{equation}
From this we can again prove \eqref{E:Y0}.

The actual argument is substantially more complicated, because in order to verify \eqref{E:nilkey} we have to prove a variant of \eqref{E:pq0} involving variable polynomial sequences, in which case the required invariance property, analogous to \eqref{E:invariance}, is much harder to identify and prove. However, the previous summary accurately summarizes the skeleton of the proof of Theorem~\ref{T:FH}.

\subsection{Vanishing of correlations}
An important component in the proof of Theorem~\ref{T:MainMulti1}, which  covers the case of Pythagorean pairs $(m^2\pm n^2, 2mn)$, is that for aperiodic multiplicative functions several correlations vanish, a result obtained in \cite[Theorem~2.6]{FH17}.
\begin{theorem}\label{T:CorrLinear}
	Let $f\colon \N\to \U$ be an aperiodic  multiplicative function. Then
	\begin{equation}\label{E:PythVanish}
	\lim_{N\to \infty}\E_{m,n\in[N]}\,  f(m^2-n^2)\cdot \overline{f(mn)}=0.
	\end{equation}
\end{theorem}
To prove this, one first uses some standard maneuvering to bound the left side  by  a multiple of $\norm{f}_{U^3(\Z)}$, which  vanishes by  Theorem~\ref{T:FH}. A similar argument leads to the following  more general result: If $f_1\colon \N\to \U$ is an aperiodic multiplicative function and $f_2,\ldots, f_\ell\colon \N\to \U$ are arbitrary sequences, then
$$
\lim_{N\to \infty}\E_{m,n\in[N]}\, f_1(L_1(m,n))\cdots f_\ell(L_\ell(m,n))=0
$$
for all linear forms $L_1,\ldots, L_\ell\colon \N^2\to \N$ such that $L_1$ is not a multiple of $L_2,\ldots, L_\ell$.

To handle Pythagorean pairs of the form $(2mn, m^2+n^2)$, we need the following result:
\begin{theorem}\label{T:CorrSumSquares}
	Let $f\colon \N\to \U$ be an aperiodic completely  multiplicative function. Then
	\begin{equation}\label{E:sumssquares}
	\lim_{N\to \infty}\E_{m,n\in[N]}\,  f(mn)\cdot \overline{f(m^2+n^2)}=0.
	\end{equation}
\end{theorem}
To prove this, we first represent    the sequence $(m,n)\mapsto f(m^2+n^2)$ as a multiplicative function on the Gaussian integers. To do this, we let $g\colon \Z[i]\to \U$ by $g(z):=\overline{f(|z|^2)}$, and note that for $z=m+in$ we have $g(z)=\overline{f(m^2+n^2)}$.
Then \eqref{E:sumssquares}  takes the following form
$$
\lim_{N\to \infty} \E_{m,n\in [N]} \, f(\Re(z)) \cdot f(\Im(z))\cdot g(z)=0.
$$
Using a variant of Theorem~\ref{T:DK} for bounded sequences on the Gaussian integers, it turns out that the last identity follows if we show that
$$
\lim_{N\to \infty}\E_{z\in R_N}\,  f(\Re(pz)) \cdot f(\Im(pz))\cdot
\overline{f(\Re(qz))} \cdot \overline{f(\Im(qz))}=0
$$
holds for all but finitely many distinct  non-real Gaussian primes $p,q$, where
$$
R_N:=\{z=m+in\colon m,n\in[N]\}.
$$
Using some elementary reasoning on the ring $\Z[i]$, it is not hard to deduce  that this  follows from
$$
\lim_{N\to \infty}\E_{m,n\in[N]} \,  f(L_1(m,n)) \cdot f(L_2(m,n))\cdot \overline{f(L_3(m,n))} \cdot \overline{f(L_4(m,n))}=0,
$$
where $f$ is aperiodic   and $L_1,L_2,L_3,L_4\colon \N^2\to \N$ are pairwise independent linear forms (recall that  we extend $f$ to the negative integers as an even function if needed).
Using standard arguments, this follows from the $U^3$-uniformity of $f$, i.e. the $s=3$ case of Theorem~\ref{T:FH}.
	
	\subsection{More general binary quadratic forms}
The results of the previous subsection can be extended to a more general setting, in which one considers the correlations of an aperiodic completely multiplicative function involving several linear forms and at most one irreducible quadratic form. A good representative example is the identity
	\begin{equation}\label{E:sumssquaresgeneral}
	\lim_{N\to \infty}\E_{m,n\in[N]} \, f(mn)\cdot \overline{f(m^2+dn^2)}=0,
\end{equation}
where $d\in \Z$ is not a square of an integer.
The case in which the expression $m^2+dn^2$ is replaced by an arbitrary irreducible binary quadratic form  $P(m,n)$, is only marginally more complicated.

Suppose first that   $d>0$. If $d\equiv 1, 2 \pmod{4}$,  then   $m^2+dn^2$
is a norm form on the ring of integers $\Z[\sqrt{-d}]$, which is a unique factorization domain. The argument employed in the case of $d = 1$ can be applied without significant alterations to this case, even in the context of a non-principal ideal domain. If $d\equiv 3 \pmod{4}$, then $\Z[\sqrt{-d}]$ is not a  unique factorization domain  and this causes problems for our argument.
When, say $d=3$, this issue can be circumvented  by working with the binary form $m^2+mn+n^2$ instead of the form $m^2+3n^2$. The advantage now is that this is a norm form on the ring $\Z[(1+\sqrt{-3})/2]$, which is a   unique factorization domain. The identity
$$
m^2+3n^2=(m-n)^2+(m-n)(2n)+(2n)^2,
$$
and a ``change of variables'' argument  are then employed to transfer results from one binary form to the other.

The  case $d<0$ presents a more substantial challenge. Let us assume that $d=-2$ and explain what is the issue in this case.
An element $z\in \Z[\sqrt{2}]$  is said to be
 {\em $C$-regular} if, for all $N\in \N$, the following condition is satisfied
$$
z^{-1}  R_{N}\subset R_{C|\mathcal{N}(z)|^{-\frac{1}{2}}N},
$$
where
 $R_N:=\{m+n\sqrt{2}\colon  m,n\in[N]\}$. In order to prove a useful variant of Theorem~\ref{T:DK},  we would like to know that for some fixed $C>0$, all elements of the ring $\Z[\sqrt{2}]$
 are $C$-regular. Unfortunately, this property does not hold, for the reason that
 $\Z[\sqrt{2}]$ has infinitely many units, so every element has infinitely many associates, and not all of them have good regularity properties. What allows us to overcome this obstacle is a result proved by Sun in \cite{Su23}, namely that for a fixed $C>0$, every element has a $C$-regular associate. This enables  us to prove a variant of  Theorem~\ref{T:DK}
 for $\Z[\sqrt{2}]$ and more general quadratic integer rings, which is applicable for our purposes.

 By combining the previous ideas we prove the following result in \cite[Theorem~3.1]{FKM24}:
 \begin{theorem}\label{T:CorrGeneral}
	Let $s\in \N$ and $f_1,\ldots, f_s,g\colon \N\to \S^1$  be  completely multiplicative functions   that are  extended to even sequences in $\Z$ and suppose  that $f_1$ is aperiodic.  Let also
$L_1,\ldots, L_s$ be linear forms in two variables with integer coefficients  and suppose that either $s=1$ and $L_1$ is non-trivial, or $s\geq 2$ and the forms $L_1,L_j$ are linearly independent for $j=2,\ldots,s$.  Finally, suppose that the form
$$
P(m,n):=\alpha m^2+\beta mn+\gamma n^2,
$$
where $\alpha,\beta, \gamma\in \Z$,   	  	   is
irreducible. Then
$$
\lim_{N\to\infty} \E_{m,n\in [N]}\, \prod_{j=1}^sf_j(L_j(m,n))\cdot g(P(m,n))=0.
$$
 \end{theorem}
For the partition regularity applications we have in mind, we  only use the case $s=2$.

	\section{Pretentious multiplicative functions - Concentration estimates}\label{S:concentration}
Now that we have at our disposal the vanishing property  of the correlations of aperiodic multiplicative functions given by Theorem~\ref{T:CorrGeneral}, we turn our attention to the study of those that are not aperiodic. In this case, we prove an approximate periodicity property that is analytically expressed in the form of the concentration estimates given in Theorems~\ref{T:concestlinear}, \ref{T:concquadratic}, and \ref{T:concquadraticgeneral}. In the non-aperiodic case, this will allow us to substantially simplify the averages that appear in  Theorems~\ref{T:MainMulti1} and \ref{T:MainMulti2}.

	\subsection{Pretentious  distance and characterization of  aperiodicity}\label{SS:pretentious}
	Our first goal is to give a useful characterization of non-aperiodic multiplicative functions.
  The following are two potential obstructions to aperiodicity:
	\begin{enumerate}
		\item (Dirichlet characters)
	A {\em Dirichlet character} is  a periodic completely multiplicative function. Then  $\chi(Qn+1)=1$ whenever $Q$ is a multiple of a period of $\chi$ and \eqref{E:aperiodic} fails  for $a=Q$ and $b=1$.
	
	\item (Archimedean characters) An {\em Archimedean character} is a completely multiplicative function of the form  $n\mapsto n^{it}$, $n\in \N$, where $t\in \N$. Then
	$\E_{n\in [N]} \, n^{it}$ is asymptotically     equal to $N^{it}/(1+it)$ and so
	\eqref{E:aperiodic} fails  for $a=1$ and $b=0$.
		\end{enumerate}
	It is a surprising and very useful  fact that   these are the only obstructions to aperiodicity. A completely multiplicative function that is not aperiodic is, in a sense that will be made precise in Theorem~\ref{T:DH} below, close to
	a product of a Dirichlet character and an Archimedean character.

	Following Granville and Soundararajan~\cite{GS07,GS08,GS23}, we define the notion of pretentiousness
	and a related distance between multiplicative functions.
	\begin{definition}\label{D:Pretentious}
		If  $f,g\colon \N\to \U$ are multiplicative functions, we define their distance  as
		\begin{equation}\label{E:D}
			\D^2(f,g):=\sum_{p\in \P} \frac{1}{p} \big(1-\Re(f(p)\cdot \overline{g(p)})\big).
		\end{equation}
		(Note that for $f,g\colon \N\to  \S^1$ we have
		$
		\D^2(f,g)=\frac{1}{2}\sum_{p\in \P} \frac{1}{p} |f(p)- g(p)|^2.
		$)
	\end{definition}	
	It can be shown (see \cite{GS08} or \cite[Section~2.1.1]{GS23}) that $\D$ satisfies the triangle inequality
	$$
	\D(f, g) \leq \D(f, h) + \D(h, g)
	$$
	for all  $f,g,h\colon \N\to \U$.
	Also, for all  $f_1, f_2, g_1, g_2\colon \N\to \U$,  we have (see
	\cite[Lemma~3.1]{GS07})
	\begin{equation}\label{E:Df1f2}
		\D(f_1f_2, g_1g_2) \leq \D(f_1, g_1) + \D(f_2, g_2).
	\end{equation}
	
	We say that  {\em $f$ pretends to be $g$} and write $f\sim g$  if  $\D(f,g)<+\infty$.
		It follows from \eqref{E:Df1f2} that if $f_1\sim g_1$ and $f_2\sim g_2$, then $f_1f_2\sim g_1g_2$.

	The next result is a direct consequence of the characterization of aperiodicity given by  \eqref{E:AperiodicChi} and the
	Wirsing-Hal\'asz mean value theorem \cite{Hal68}, which  gives necessary and sufficient conditions for a multiplicative function to have mean value  $0$.
	\begin{theorem}\label{T:DH}
		A  completely multiplicative function  $f\colon \N\to \U$ is  not aperiodic if and only if   there exists
		 $t\in \R$ and Dirichlet character $\chi$ such that  $f\sim \chi\cdot n^{it}$, i.e.
		$
		 \D(f,\chi\cdot n^{it})<+\infty.
		$
	\end{theorem}
 This motivates  the following definition:
 \begin{definition}
 		We say that $f$ is {\em pretentious,} if $f\sim  \chi \cdot n^{it}$ for some $t\in \R$ and Dirichlet character $\chi$.
 \end{definition}
	The value of  $t$ is uniquely determined; this follows
	from \eqref{E:Df1f2} and the fact that $n^{it}\not\sim \chi$
	for every non-zero $t\in \R$  and  Dirichlet character $\chi$ (see for example \cite[Corollary~11.4]{GS23}  or \cite[Proposition~7]{GS08}).
	
	It follows from Theorem~\ref{T:DH} that if a completely multiplicative function $f\colon \N\to \U$ is not aperiodic, then it is pretentious.

	\subsection{Examples of pretentious multiplicative functions}\label{SS:pretex}
We give several illustrative examples of pretentious multiplicative functions in order to
cover all the possible varying behaviors that are relevant to our study.
In the  examples that follow, we describe the multiplicative function only on prime numbers and assume that it is extended to a completely
multiplicative function on the positive integers.
\begin{enumerate}
		\item\label{i} (Modified Dirichlet characters) If $\chi$ is a  Dirichlet character, then $\chi(p)=0$ if and only if $p$ divides the minimal period of $\chi$ and $\chi(p)$ is a root of unity otherwise.  We let
	$\tilde{\chi}\colon \N\to \S^1$
	 be the completely multiplicative function satisfying
	$$
	\tilde{\chi}(p):=
	\begin{cases}
		\chi(p), &\quad \text{if } \chi(p)\neq 0\\
		1, &\quad \text{if } \chi(p)=0.
	\end{cases}
	$$
	Then $\tilde{\chi}$ is not necessarily periodic and $\tilde{\chi}\sim \chi$.

		\item  \label{ii}  (Finite support) Let $f(p_1)=\cdots =f(p_\ell):=-1$ and $f(p)=1$ for $p\notin \{p_1,\ldots, p_\ell\}$. Then $f\sim 1$.
	
	\smallskip

\item  \label{iii}  (Infinite support $1$-pretentious) Let $f(p):=-1$  for $p\in \P_0$ and $f(p)=1$ for $p\notin \P_0$ where $\P_0\subset \P$ is  infinite and
$\sum_{p\in \P_0}\frac{1}{p}<+\infty$. Then $f\sim 1$.

	\smallskip
	
	\item \label{iv} ($1$-pretentious oscillatory)   $f(p):=e(1/ \log\log{p})$, $p\in \P$. Then $f\sim 1$ but  it turns out that $f$ does not have a mean value, in fact, $\E_{n\in [N]}\, f(n)$ is asymptotically equal to $e(w(N))$ where
	$w(N)=c\log\log\log N$ for some $c>0$.
	
	\item \label{v} (Archimedean factor)  $f(n):=n^{it}\cdot g(n)$ for some $t\in \R$ and $g$  as in the previous examples.  Then $f\sim \chi \cdot n^{it}$ if $g$ is as in \eqref{i} and $f\sim  n^{it}$ if $g$ is as in \eqref{ii}-\eqref{iv}.
\end{enumerate}

We mention some   features of the previous examples that will motivate the concentration results to be proved in the next subsections.

In the examples~\eqref{i} and \eqref{ii}, there exists $Q_0$ such that if $Q\mid Q_0$, then
$$
f(Qn+1)=1 \quad \text{for every } n\in \N.
$$
 Indeed, take $Q_0$ be the period of $\chi$ in the first case and $Q_0:=p_1\cdots p_\ell$ in the second.
	
	In example~\eqref{iii}, let $\varepsilon>0$ and $K=K(\varepsilon)$ be such that
	$\sum_{p\in \P_0, \, p> K}\frac{1}{p}<\varepsilon$.  For $Q_0:=\prod_{p\leq K} p$ we have that if $Q_0\mid Q$, then  $\overline{d}(n\in \N\colon f(Qn+1)\neq 1)\leq \varepsilon$ and hence
	$$
	\limsup_{N\to\infty}\E_{n\in [N]}|f(Qn+1)-1|\ll \varepsilon.
	$$

	In example~\eqref{iv}, let $\varepsilon>0$ and $K:=K(\varepsilon)$ be such that
$\sum_{p\in \P, \, p> K} \frac{1}{p}\big(1-\Re(f(p))\big)<\varepsilon$.  For $Q_0=\prod_{p\leq K} p$ we have that if $Q_0\mid Q$, then
$$
\limsup_{N\to\infty}\E_{n\in [N]}|f(Qn+1)-e(w(N))|\ll \varepsilon,
$$
where $w(N)$ increases to infinity slowly, in fact like $c\log\log\log{N}$ for some $c>0$.

	In example~\eqref{v}, we have
$$
\lim_{N\to\infty}\E_{n\in [N]}|f(Qn+1)-(Qn)^{it}\cdot g(Qn+1)|=0.
$$
Combining this  with the previous estimates for $g$ we get analogous results for $f$.

The general principle that can be derived from these examples is that if $f\sim \chi$, then  for highly divisible values of $Q$ the values  $f(Qn+1)$ concentrate around  $1$ or  a slowly growing oscillatory factor.
The next subsection makes this principle explicit. Later, we also give variants of it that apply to arbitrary binary quadratic forms in place of $n$.

\subsection{Linear concentration estimates - Statements}\label{SS:conclin}	
To get a sense of the concentration estimates we will be using in a simpler setting,  let us first consider the case  where a pretentious multiplicative function takes on finitely many values on the primes.
\begin{theorem}\label{T:fv}
Let $f\colon \N\to \U$ be a pretentious multiplicative function and suppose that the set $\{f(p)\colon p\in \P\}$ is finite. Then for every  $\varepsilon>0$ there exists $Q_0=Q_0(\varepsilon)$ such that if $Q\in \N$ is such that $Q_0\mid Q$, then
	\begin{equation}\label{E:QN1}
	\limsup_{N\to\infty}\E_{n\in[N]}|f(Qn+1)-1|\leq \varepsilon.
	\end{equation}
\end{theorem}
Such concentration  results are distinctive properties of pretentious multiplicative functions. Indeed,  for every $Q\in \N$ an aperiodic multiplicative function has average $0$ on the arithmetic progression $Q\N+1$, hence \eqref{E:QN1} fails if  $\varepsilon <1$.

We also note that a related property of ``Besicovitch rational almost periodicity'', obtained in \cite[Theorem~6]{DD82} for pretentious multiplicative functions $f\colon \N\to \U$ that take finitely many values on primes, states that
for every $\varepsilon>0$ there exists a periodic sequence $a_\varepsilon$ such that
$$
\limsup_{N\to\infty} \E_{n\in[N]}|f(n)-a_\varepsilon(n)|\leq \varepsilon.
$$
Unfortunately, for large values of $Q$, this estimate alone cannot give us useful information about the values of $f(Qn+1)$, and so it is not helpful for our applications.

For our purposes we need an extension  of Theorem~\ref{T:fv},  that covers all pretentious multiplicative functions
and also allows for some uniformity in our choice of $Q$. To state it, we have to introduce some notation. For $K\in \N$ let
\begin{equation}\label{E:PhiK}
	\Phi_K	:=\Big\{\prod_{p\leq K}p^{a_p}\colon K< a_p\leq 2K\Big\}.
\end{equation}
If $f\sim \chi\cdot n^{it}$ for some Dirichlet character $\chi$ and $t\in \R$ and
 $K,N\in \N$, we let
	 \begin{equation}\label{E:FNQ}
		F_N(f,K):=\sum_{K< p\leq N} \frac{1}{p}\,\big(f(p)\cdot \overline{\chi(p)}\cdot p^{-it} -1\big).
	\end{equation},

Since $f\sim \chi\cdot n^{it}$, we have that the limit $\lim_{N\to\infty} \Re( F_N(F,1))$
exists. But   there are cases  where  $\lim_{N\to\infty} |\text{Im} (F_N(f,1))| =+\infty$, this holds in example~\eqref{iv} of Section~\ref{SS:pretex}.

\begin{theorem}\cite[Lemma 2.5]{KMPT21}\label{T:concestlinear}
	Let $f\colon \N\to \U$ be a multiplicative function and suppose that $f\sim \chi\cdot n^{it}$ for some Dirichlet character $\chi$ and $t\in \R$.
	 Then	
	 \begin{equation}\label{E:fconc}
	 \lim_{K\to\infty} 	\limsup_{N\to\infty} \max_{Q\in \Phi_K} \E_{n\in[N]}\big|f(Qn+1)-  (Qn)^{it}\cdot \exp\big(F_N(f,K)\big)\big|=0,
	 \end{equation}
	where $\Phi_K$ and $F_N(f,K)$ are as in \eqref{E:PhiK} and  \eqref{E:FNQ} respectively.
\end{theorem}
Note that in the case of multiplicative functions  that take finitely many values on primes, we have $t=0$ and
the series defining $F_N(f,K)$ can be shown to converge, so in our statement we can replace
the term $F_N(f,K)$ with $0$. This gives a statement similar to the one in Theorem~\ref{T:fv}.

\subsection{Linear concentration estimates - Proof} There are two main steps in the proof of Theorem~\ref{T:concestlinear}. In the first step we deal with additive functions on the complex numbers and prove the concentration estimate \eqref{E:hconc} under the assumption \eqref{E:hpre}. We  describe the key steps of this argument, the reader will find the missing details, for example, in \cite[Chapter III.3.2]{T95}. In the second step, given a pretentious multiplicative function $f$,  we apply the previous results for the additive  function $h$ defined in \eqref{E:hf} to derive the required concentration estimate for $f$.

\medskip

 {\bf Step 1 (Tur\'an-Kubilius inequality).}  The function  $h\colon \N\to \C$ is {\em additive} if it satisfies
$$
h(mn)=h(m)+h(n)\quad  \text{whenever } (m,n)=1.
$$
For simplicity we also assume that
\begin{equation}\label{E:2zero}
h(p^k)=0 \quad \text{for every } p\in \P, \,  k\geq 2,
\end{equation}
and remark that this simplifying   assumption can be easily removed.
The Tur\'an-Kubilius  inequality implies that if
\begin{equation}\label{E:hpre}
\sum_{p\in \P}\frac{|h(p)|^2}{p}<\infty,
\end{equation}
then
\begin{equation}\label{E:hconc}
\lim_{K\to\infty} 	\limsup_{N\to\infty} \max_{Q\in \Phi_K} \E_{n\in[N]}|h(Qn+1)-  H_N(h,K)|^2=0,
\end{equation}
where  $\Phi_K$ is as in \eqref{E:PhiK} and
$$
H_N(h,K):=\sum_{K<p\leq N}\frac{h(p)}{p}.
$$
(As it will turn out, for $Q\in \Phi_K$, $H_N(h,K)$ is approximately equal to  $\E_{n\in [N]}\, h(Qn+1)$).

 We mention here the important ingredients of the proof of \eqref{E:hconc}.  In the next subsections we list the essential changes that need to be made to get analogous extensions for  binary quadratic forms.

 Let
 $$
 w_{Q,N}(p, q):=\frac{1}{N}\sum_{n\in [N], \, p,q\mid\mid Qn+1} 1,
 $$
 where  $p\mid\mid Qn+1$ means that $p\mid Qn+1$ but  $p^2\nmid Qn+1$.
To prove \eqref{E:hconc},   we expand the square  and  compute it using the identity
$$
h(Qn+1)=\sum_{p\mid\mid Qn+1} h(p)=\sum_{p>K, \, p\mid\mid Qn+1} h(p),
$$
which follows from the additivity of $h$, our standing simplifying assumption $h(p^k)=0$ for $k\geq 2$,  and the fact that if $Q\in \Phi_K$ we have  $p\nmid Qn+1$ for all primes $p\leq K$.
 We deduce  that
\begin{equation}\label{E:h1}
\E_{n\in[N]}\, h(Qn+1)= \sum_{K<p\leq N}  w_{Q,N}(p,p)\cdot h(p)
\end{equation}
and
\begin{equation}\label{E:h2}
\E_{n\in[N]} |h(Qn+1)|^2= \sum_{K<p, q \leq N}  w_{Q,N}(p,q)\cdot h(p)\cdot \overline{h(q)}.
\end{equation}

By direct computation, we get   the approximate identities
$$
w_{Q,N}(p, p)= \frac{1}{p}\Big(1-\frac{1}{p}\Big) +O\Big(\frac{1}{N}\Big),
$$
and
$$
w_{Q,N}(p, q)= \frac{1}{pq}\Big(1-\frac{1}{p}\Big)\Big(1-\frac{1}{q}\Big)+O\Big(\frac{1}{N}\Big),
$$
which hold for all primes $p,q$ with $p\neq q$. Importantly, for $p\neq q$ we
have
\begin{equation}\label{E:w3}
w_{Q,N}(p, q)=w_{Q,N}(p, p)\cdot w_{Q,N}(q, q) +O\Big(\frac{1}{N}\Big).
\end{equation}
(The fact that the error term  in \eqref{E:w3} does not depend on $Q$ is  responsible for the uniformity with respect to $Q$ in \eqref{E:fconc}.)
Using  \eqref{E:h1}, \eqref{E:h2}, \eqref{E:w3}, it is possible to show that \eqref{E:hconc} holds. To do this, we use that terms of the form $ \sum_{K<p\leq N}\frac{|h(p)|^2}{p}$ have negligible contribution in \eqref{E:hconc}  because of  \eqref{E:hpre}. Also error terms of the form  $O\big(\frac{1}{N}\big)$
contribute a total error
$
O\big(\sum_{p,q\colon pq\leq QN+1}\frac{1}{N}\big)=O(Q \log\log{N}/\log{N}),
$ which is also negligible for  \eqref{E:hconc}.
This last fact is no longer true for polynomials in two variables and respective averages, because there are too many error terms, which causes significant problems in our treatment of analogous concentration estimates for sums of two squares and more general irreducible binary quadratic forms.

\medskip

 {\bf Step 2 (From additive to multiplicative).}  The second step is to pass the previous  result from ``pretentious''  additive functions to pretentious multiplicative functions.
To do this, we first reduce to the case where $f\sim 1$ and  then define the additive function
$h\colon \N\to \C$ on prime powers as follows
\begin{equation}\label{E:hf}
h(p^k):=f(p^k)-1, \quad p\in \P, \, k\in \N.
\end{equation}
We also make the simplifying assumption $f(p^k)=1$ for $k\geq 2$, hence
\eqref{E:2zero} holds.  We can remove this   assumption  using that $f\sim 1$, which implies  that the contribution of prime powers greater than $1$ is negligible in \eqref{E:fconc}.

Since $f\sim 1$, we get that \eqref{E:hpre} holds and hence  \eqref{E:hconc} holds from  Step 1.
Our goal is to relate the values of $f(Qn+1)$ with those of $h(Qn+1)$ and derive \eqref{E:fconc} form \eqref{E:hconc}. We  will use the approximate identity
\begin{equation}\label{E:zez}
z=e^{z-1}+O(|z-1|^2),
\end{equation} which holds for $|z|\leq 1$. Note that since $f\sim 1$ we have that $f(p)-1$ is close to $0$ for ``most'' primes $p$, so  for $z_p:=f(p)$  we have that   \eqref{E:zez} gives a sensible approximation.
We deduce  that
$$
f(Qn+1)= \prod_{p\mid\mid  Qn+1, p>K} f(p)=\prod_{p\mid\mid Qn+1, p>K}\big(\exp(h(p))+O(|h(p)|^2)\big),
$$
where we used that  for $Q\in \Phi_K$ only  primes greater than $K$ divide $Qn+1$.  If we ignore the
error terms on the right side  (this can be done using \eqref{E:hpre}), we get that it is equal to
$$
\prod_{p\mid\mid Qn+1}\exp(h(p))=\exp(h(Qn+1)),
$$
where we used \eqref{E:2zero}.
Summarizing, we have just shown  that in establishing \eqref{E:fconc}, we can replace $f(Qn+1)$ with
$\exp(h(Qn+1))$.
So to prove \eqref{E:fconc}, it remains to show that (note that since $f\sim 1$ we have $t=0$)
 $$
\lim_{K\to\infty} 	\limsup_{N\to\infty} \max_{Q\in \Phi_K} \E_{n\in[N]}\big|  \exp(h(Qn+1))-\exp\big(F_N(f,K)\big)\big|=0.
$$
Using the estimate $|e^z-e^w|\ll |z-w|$,  which holds for $\Re(z), \Re(w)\leq 0$, and the Cauchy-Schwarz inequality,  it
suffices to show that
$$
\lim_{K\to\infty} 	\limsup_{N\to\infty} \max_{Q\in \Phi_K} \E_{n\in[N]}|  h(Qn+1)-F_N(f,K)|^2=0.
$$
Since  $F_N(f,K)=H_N(h,K)$,   this follows from \eqref{E:hconc}.

	\subsection{Concentration estimates along sums of two squares}\label{SS:concquad}
	To prove that the pair of polynomials $2mn$, $m^2+n^2$ is density regular (which implies that $x^2+y^2=z^2$ is partition regular with respect to $y,z$), we need to verify
 Theorem~\ref{T:MainMulti1} when $\alpha=\beta=1, \gamma=2$.
To do this,  we need
   a variant of Theorem~\ref{T:concestlinear} that deals with concentration estimates along sums of two squares. We give the precise statement next.

	If $f\sim \chi\cdot n^{it}$ for some Dirichlet character $\chi$ and $t\in \R$, and
	$K,N\in \N$, we let
	\begin{equation}\label{E:FNQsquares}
		G_N(f,K):=\sum_{K< p\leq N, \,  p\equiv 1 \! \! \! \!\!  \pmod{4}} \frac{1}{p}\,\big(f(p)\cdot \overline{\chi(p)}\cdot p^{-it} -1\big).
	\end{equation}
	This expression is similar to  $F_N(f,K)$  in \eqref{E:FNQ}, but we omit   primes that are congruent to $3 \mod{4}$ in the summation.
 This is due to the fact that no prime of this form is a sum of two squares, which ultimately makes the contribution of such primes negligible for our purposes.
	\begin{theorem}\cite[Proposition~2.11]{FKM23}\label{T:concquadratic}
		Let $f\colon \N\to \U$ be a multiplicative function and suppose that $f\sim \chi\cdot n^{it}$ for some Dirichlet character $\chi$ and $t\in \R$.
		Then	
		\begin{multline}\label{E:fconc2}
			\lim_{K\to\infty} 	\limsup_{N\to\infty} \max_{Q\in \Phi_K} \E_{m,n\in[N]}\, \big|f\big((Qm+1)^2+(Qn)^2\big)-\\  Q^{2it}\cdot (m^2+n^2)^{it}\cdot \exp\big(G_N(f,K)\big)\big|=0,
		\end{multline}
		where $\Phi_K$ and $G_N(f,K)$ are as in \eqref{E:PhiK} and  \eqref{E:FNQsquares} respectively.
	\end{theorem}
If $f$ is pretentious and  takes finitely many values in the primes,  then  it can be shown that $f\sim \chi$ for some Dirichlet character $\chi$ and the series defining $G_N(f,K)$ converges as $N\to\infty$, so the term
	 $ Q^{2it}\cdot (m^2+n^2)^{it}\cdot \exp\big(G_N(f,K))$ in \eqref{E:fconc2} can be replaced by  $1$.
	
Our strategy is similar to the one used to prove the linear concentration estimates in Section~\ref{SS:conclin}. However, to carry out the first step, which is the Tur\'an-Kubilius estimate for additive functions taken along sums of two squares,  there are some additional important difficulties that we must resolve.

Our goal is to show that  if
\begin{equation}\label{E:hprsquares}
	\sum_{p\in \P}\frac{|h(p)|^2}{p}<\infty,
\end{equation}
and $h_N$ is the restriction of $h$ on primes less than $N$, i.e.
$$
h_N(p^k):= \begin{cases}h(p^k), \quad  &\text{if }\     p\leq N\\
	0, \quad &\text{otherwise}  \end{cases},
$$
then  we have the following variant of the Tur\'an-Kubilius inequality \eqref{E:hconc} that deals with sums of two squares
\begin{equation}\label{E:hconsquares}
	\lim_{K\to\infty} 	\limsup_{N\to\infty} \max_{Q\in \Phi_K} \E_{m,n\in[N]}|h_N((Qm+1)^2+(Qn)^2)-  H_N(h,K)|^2=0,
\end{equation}
where
$$
H_N(h,K):=\sum_{K<p\leq N, \, p\equiv 1  \! \! \! \!\!  \pmod{4} }\frac{h(p)}{p}.
$$
We could replace  $h_N$ with $h$  in \eqref{E:hconsquares},  but  because of some error terms that appear in the next step  of the argument,  such a result  would not be more useful.
For convenience, we also assume  that $h(p)=0$ for all primes that are congruent to $3 \pmod{4}$ and $h(p^k)=0$ for all primes $p$ and $k\geq 2$; both conditions can be easily removed.

Following the plan used to prove the  linear Tur\'an-Kubilius inequality  \eqref{E:hconc}, we  must  compute
$$
w_{Q,N}(p, q):=\frac{1}{N^2}\sum_{m,n\in[N], \, p,q\mid\mid (Qm+1)^2+(Qn)^2} 1
$$
whenever $p,q \equiv 1  \pmod{4}$.
Using  that $-1$ is a quadratic residue $\!  \! \mod{p}$ whenever $p \equiv 1 \pmod{4}$ is a quadratic residue and the Chinese remainder theorem, we find that
	$$
	w_{Q,N}(p, p)= \frac{2}{p}\Big(1-\frac{1}{p}\Big)^2 +O\Big(\frac{1}{N}\Big),
	$$
	and
	\begin{equation}\label{E:pq}
	w_{Q,N}(p, q)= \frac{4}{pq}\Big(1-\frac{1}{p}\Big)^2\Big(1-\frac{1}{q}\Big)^2+O\Big(\frac{1}{N}\Big),
\end{equation}
	which holds for all primes $p,q$ with $p\neq q$ that are congruent to $1\pmod{4}$.  Importantly, for primes $p\neq q$,  congruent to $1\pmod{4}$,  we
again 	have the approximate identity
	\begin{equation}\label{E:w3'}
		w_{Q,N}(p, q)=w_{Q,N}(p, p)\cdot w_{Q,N}(q, q) +O\Big(\frac{1}{N}\Big).
	\end{equation}
In contrast to the linear case, in our case, because we are averaging over two variables $m,n$,  the errors incurred by terms of the form  $O\big(\frac{1}{N}\big)$ are too many and cannot be easily neglected. To deal with this problem, we decompose $h$ as a sum $h_1+h_2$ of two additive functions, ``supported'' on
 primes less than $\sqrt{N}$ and primes greater than $\sqrt{N}$, respectively. The $O\big(\frac{1}{N}\big)$ errors turn out to be   negligible for $h_1$. For $h_2$ they  can be handled  using the following simple but important fact:
\begin{lemma}\label{L:keypq}
If $p,q$ are distinct primes congruent to  $1\pmod{4}$, then
$$
w_{Q,N}(p,p)\ll \frac{Q^2}{p}, \quad w_{Q,N}(p,q)\ll \frac{Q^2}{pq},
$$
where the implicit constant is absolute.
\end{lemma}
The proof of this fact is simple and is based on the  elementary estimate
$$
\sum_{k\leq n}r_2(k)\ll n,
$$
where $r_2(n)$  is the number of representations of $k$ as a sum of two squares.
Lemma~\ref{L:keypq}  gives better upper bounds than \eqref{E:pq} when $pq$ is much larger than $N$,  bounds that turn out to be sufficient to show that the contribution of $h_2$ in \eqref{E:fconc2} is negligible
(\eqref{E:hprsquares} with $h_2$ instead  of $h$ is also used here).

The second step of the proof is the transition from additive to multiplicative functions. We follow the same approach as in the linear case described above. We first reduce to the case $f\sim 1$, in which case we get that \eqref{E:fconc2} holds with $f_N$ instead of $f$, where
$f_N$ is the restriction of $f$ on primes less than $N$, i.e.
	$$
	f_N(p^k):= \begin{cases}f(p^k), \quad  &\text{if }\     p\leq N\\
		1,  \quad &\text{otherwise}  \end{cases}.
	$$
The  error terms incurred in this transition  allow us to deal only with additive and multiplicative functions supported  on primes $p\leq N$, this is the reason why in \eqref{E:hconsquares} a larger range of primes in the support of $h_N$ would not have been helpful.

Finally, to replace  $f_N$ with  $f$ in \eqref{E:fconc2}, we use   the factorization	 $f=f_N\cdot g_N$, where $g_N$ is the restriction of $f$ on primes greater than $N$, defined in a similar way to $f_N$. 	We then use  a slight variant of   Lemma~\ref{L:keypq}  together with our hypothesis $f\sim 1$ to show that the contribution of $g_N$ in \eqref{E:fconc2} is negligible and complete the proof of Theorem~\ref{T:concquadratic}.

\subsection{Concentration estimates for general binary quadratic forms}\label{SS:conentrationgeneral}
To prove Theorems~\ref{T:MainMulti1} and \ref{T:MainMulti2},
which are needed to cover partition and density regularity results of the more general quadratic equations $ax^2+by^2=cz^2$,
we need to establish variants of Theorem~\ref{T:concquadratic} in which the  place of $m^2+n^2$ is taken by the more general irreducible binary quadratic form $P(m,n)=\alpha m^2+\beta mn +\gamma n^2$, where $\alpha,\beta,\gamma\in \Z$.
To simplify the notation a bit, we give the exact statement when $P(m,n)=m^2+dn^2$ and $d\in \Z$ is a square-free integer, the general case is similar and the exact statements can be found in \cite[Proposition~4.9]{FKM24}.
 We first have to modify the definition of $G_N(f,K)$ from \eqref{E:FNQsquares}.
If  $\legendre{a}{p}$ denotes the Legendre symbol and  $d\in \Z$ is square-free, we let
 	\begin{equation}\label{E:Pd}
 	\CP_{d}:=\Big\{p\in \P\colon  \legendre{-d}{p}=1 \Big\}.	
 \end{equation}

 For example,
 $\CP_1:=\{ p\equiv 1\!  \pmod{4}\}$,
 $\CP_2:=\{p\equiv 1, 3\!\pmod{8}\}$,
 $\CP_{-2}:=\{ p\equiv 1, 7\!\!\pmod{8}\}$.
 It can be shown that  for square-free $d$ the set $\CP_{d}$  has relative density $1/2$ on the primes and it always consists of the set of  primes that belong to  finitely many congruence classes. For example, using quadratic reciprocity, we get that  if $q$ is an odd prime, then
 $$
 \CP_q=\{ p\equiv \pm k^2\! \! \pmod{4q}\colon  k\in [4q]
 \text{ is odd and coprime to } 4q\}.
 $$

  If $f\sim \chi\cdot n^{it}$ for some Dirichlet character $\chi$ and $t\in \R$,  and
 $K,N\in \N$, we let
 	\begin{equation}\label{E:GNd}
 	G_{d,N}(f,K):= 2\, \sum_{K< p\leq N,\,  p\in \CP_d}\, \frac{1}{p} \,(f(p)\cdot \overline{\chi(p)}\cdot n^{-it} -1).
 \end{equation}

	\begin{theorem}\cite[Proposition~4.9]{FKM24}\label{T:concquadraticgeneral}
	Let $f\colon \N\to \U$ be a multiplicative function and suppose that $f\sim \chi\cdot n^{it}$ for some Dirichlet character $\chi$ and $t\in \R$.
	Then	for every square-free $d\in \Z$ such that the form $m^2+dn^2$ is irreducible,  we have
\begin{multline*}
		\lim_{K\to\infty} 	\limsup_{N\to\infty} \max_{Q\in \Phi_K} \E_{m,n\in[N]}\big|f\big((Qm+1)^2+d\,  (Qn)^2\big)- \\ Q^{2it}\cdot (m^2+dn^2)^{it}\cdot \exp\big(G_{d,N}(f,K)\big)\big|=0,
\end{multline*}
	where $\Phi_K$ and $G_{d,N}(f,K)$ are as in \eqref{E:PhiK} and  \eqref{E:GNd} respectively.
\end{theorem}
We also prove variants of concentration estimates along more general lattices of the form $\{(Qm+a,Qn+b)\colon m,n\in \N \}$
for suitable $a,b\in \Z_+$; such results are not needed for the proof of Theorem~\ref{T:MainMulti1}, but turn out to be essential for the proof of Theorem~\ref{T:MainMulti2}.

To prove Theorem~\ref{T:concquadraticgeneral}, we follow the approach used to prove Theorem~\ref{T:concquadratic}, but there are quite a few additional difficulties.
The most serious problem is to find a variant of Lemma~\ref{L:keypq} that allows us to discount the contribution of large primes in the various steps of the proof.
For $d>0$, the  most problematic case is when  the corresponding quadratic integer ring is not a principal ideal domain, and
  additional, more substantial problems arise when  $d<0$. For example, suppose that $d<0$ is a square-free integer such that $d\equiv 1,2 \pmod{4}$ and   we want to establish for all primes $p\in \CP_{d}$ an estimate of the form
\begin{equation}\label{E:key-2}
|\{m,n\in [N]\colon p\mid m^2+dn^2\}|\ll_d \frac{N^2}{p},
\end{equation}
where the implicit constant depends only on $d$.
 A priori it is not clear that such an estimate holds, because the counting is affected by the number of the solutions of  Pell's equation $m^2+dn^2=1$  in the box $[1,N]\times [1,N]$, which grows as $N\to\infty$, roughly like $\log{N}$ (unlike the case $d>0$, where the corresponding number  is bounded by $4$).

  Fortunately,  we can  perform a careful count and use classical but non-trivial facts from  algebraic number theory to show that  \eqref{E:key-2} indeed holds. We first note that the left side in \eqref{E:key-2} is bounded by
 $$
 \sum_{r\leq \frac{2N^2}{p}} C_{d,N}(pr),
 $$
 where $C_{d,N}(k):=|\{m,n\in [N]\colon  |m^2+dn^2|=k\}|$. Using the
fact that  there exists a unit $u\in \Z[\sqrt{-d}]\subset\R$ such that
 $u>1$ and every other unit has the form $\pm u^t$ for some $t\in \Z$
  (this is the Dirichlet unit theorem for the ring $\Z[\sqrt{-d}]$), it is possible to show  that there exists  $c_d>0$ such that
 \begin{equation}\label{E:CdNk}
 C_{d,N}(k)\ll \log(c_dN/\sqrt{k})\cdot \mathcal{C}_d(k),
 \end{equation}
 where
  $$
  \mathcal{C}_d(k):=|\{ \text{ideals in } \Z[\sqrt{-d}]\text{ with norm } k\}|.
  $$
  It is important for our argument that the term $\sqrt{k}$ appears in the denominator in the right side of \eqref{E:CdNk}, our argument would not work without it. So \eqref{E:key-2} follows if we establish the bound
 \begin{equation}\label{E:estimate}
 \sum_{r\leq \frac{2N^2}{p}}\log(c_dN/\sqrt{pr})\cdot \mathcal{C}_d(pr)\ll_d \frac{N^2}{p}.
 \end{equation}
 To prove this, we crucially use the estimate
\begin{equation}\label{E:Cdpr}
 \mathcal{C}_d(pr)\leq 2\, \mathcal{C}_d(p),
 \end{equation}
  which holds for all $p\in \CP_{d}$.
This seemingly simple fact is not easy to prove from first principles, especially in the case where the corresponding ring of integers is not a principal ideal domain.
However, its proof follows easily from well-known, but non-trivial  facts in  algebraic number theory:   the multiplicativity of  $\mathcal{C}_d(k)$, the estimate
$\mathcal{C}_d(p^{k+1})\leq \mathcal{C}_d(p^k)+1$,  and the fact that
$\mathcal{C}_d(p^{k})\geq 1$ for all $p\in \CP_d$ and $k\in\N$.
The proof of \eqref{E:estimate} follows easily by using  \eqref{E:Cdpr},
 and then applying partial summation combined with  the estimate
  $$
  \sum_{r=1}^N\mathcal{C}_d(r)\ll_d N,
  $$
   which follows from the ideal counting theorem.

\section{Connecting the pieces (the $Q$-trick)}\label{S:Qtrick}
In this section we show how we can combine the ingredients of the previous two sections to prove Theorems~\ref{T:MainMulti1} and \ref{T:MainMulti2}. We start with the simplest case and then explain the necessary changes needed to cover  the increasingly more complicated cases.
\subsection{The pair $(2mn, m^2-n^2)$.} \label{SS:pair1}
We start with the  proof of Theorem~\ref{T:MainMulti1} when $\alpha=1$, $\beta=-1$, $\gamma=2$, which suffices to prove partition and density regularity of $x^2+y^2=z^2$ with respect to the variables $x,y$.
Using the positivity property \eqref{E:positivers},  it suffices to show that there exist $Q\in \N$ for which
\begin{equation}\label{E:positiveA1}
 \int_\CM 	L(f,Q)\,  d\sigma(f)>0,
\end{equation}
where
\begin{equation}\label{E:LfQ}
L(f,Q):=\lim_{N\to \infty} \E_{m,n\in[N]} \, f((Qm+1)^2-(Qn)^2)\cdot \overline{f(2(Qm+1) Qn)}.
	\end{equation}
Note that in establishing the above reduction we used the bounded convergence theorem and the non-trivial fact that  the limit defining $L(f,Q)$ exists, which follows from~\cite[Theorem~1.4]{FH16}. The proof of the existence of the limit is non-trivial, and can be avoided at the cost of making the argument  a bit more complex.

We also note that 	we could have used  in our averaging  any other lattice of the form $(Qm+a,Qn)$, where $a$ is some fixed number independent of $Q$; what is important is that on one of the two coordinates we use $Qn$  (or $Qm$) in order to later on factor out a term $f(Q)$, which is an essential element of  a key maneuver that  we call the ``$Q$-trick''.

It easily follows from  Theorem~\ref{T:CorrLinear}  that if the multiplicative function  $f$ is aperiodic, then
$$
L(f,Q)= 0 \quad \text{for every } Q\in \N.
$$
  Hence, \eqref{E:positiveA1} would follow if we show that
\begin{equation}\label{E:positiveA2}
	 \int_{\CM_p} 	L(f,Q) \,  d\sigma(f)>0 \quad \text{holds for some } Q\in \N,
\end{equation}
where
$$
\CM_p:=\{f\in \CM \text{ is pretentious} \}.
$$
Note that the integrals in \eqref{E:positiveA2} can be shown to be real.
To prove \eqref{E:positiveA2} we will use the linear concentration estimates of Theorem~\ref{T:concestlinear} to substantially simplify the value of the limit in \eqref{E:LfQ}, and then use a trick, called the $Q$-trick,  to show that for some particular choice of $Q$ we get positivity.

\subsubsection{Finitely generated case}\label{SS:fg}
To better illustrate the argument we first assume that the measure $\sigma$ is supported on the subset  of pretentious multiplicative functions that are finitely generated, i.e. the set $\{f(p)\colon p\in\P\}$ is a finite subset of the unit circle. The advantage in this case is that  $f\sim \chi$ for some Dirichlet character $\chi$ (so $n^{it}$ is no longer needed) and also the oscillatory factor $\exp\big(F_N(f,K)\big)$  in \eqref{E:fconc} can be taken to be $1$. So in this case the concentration estimate in \eqref{E:fconc} takes the much simpler form
 \begin{equation}\label{E:fconc'}
	\lim_{K\to\infty} 	\limsup_{N\to\infty} \max_{Q\in \Phi_K} \E_{n\in[N]}|f(Qn+1)-  1|=0,
\end{equation}
where $\Phi_K$  are as in \eqref{E:PhiK}. Since
$$
Qm+1\pm Qn=Q(m\pm n)+1,
$$
with a bit of fiddling (see \cite[Lemma~3.1]{FKM23}), \eqref{E:fconc'} implies that for all practical purposes the limit
$L(f,Q)$ in \eqref{E:LfQ}  can be approximated by $f(Q)\cdot A(f)$,
where
\begin{equation}\label{E:Af}
A(f):=\lim_{N\to\infty}\E_{n\in [N]} \, \overline{f(2n)}.
\end{equation}
More precisely, we have
\begin{equation}\label{E:Lfa}
	\lim_{K\to\infty}  \max_{Q\in \Phi_K}|L(f,Q)- f(Q)\cdot A(f)|=0.
\end{equation}
 It follows from this, that in order to establish \eqref{E:positiveA2} it suffices to study the integrals
 $$
 \int_{\CM_p} f(Q)\cdot A(f)\,  d\sigma(f).
 $$
A priori, however, these expressions do not satisfy any obvious positivity property, and it seems that we are at a dead end. The simple but crucial observation is that the multiplicative average of this last expression, taken over $Q\in \Phi_K$ and letting $K\to\infty$, does satisfy a nice positivity property, namely,
  \begin{equation}\label{E:id1}
\lim_{K\to\infty} \E_{Q\in \Phi_K}\int_{\CM_p} f(Q)\cdot A(f)\,  d\sigma(f)=\sigma(\{1\})>0.
 \end{equation}
 This is because   for every completely multiplicative function $f\neq 1$ we have
$$
\lim_{K\to\infty} \E_{Q\in \Phi_K}f(Q)=0,
$$
and so the bounded convergence theorem gives that \eqref{E:id1} holds.
Combining \eqref{E:Lfa} and \eqref{E:id1} we get that
$$
 \lim_{K\to\infty} \E_{Q\in \Phi_K} \int_{\CM_p} 	L(f,Q) \,  d\sigma(f)=\sigma(\{1\})>0,
$$
hence \eqref{E:positiveA2} holds.
\subsubsection{General case}  We now consider the general case of measures supported by multiplicative functions that are not necessarily finitely generated.
 In this case, suppose that  $f\sim \chi\cdot  n^{it}$ for some Dirichlet character $\chi$ and $t\in \R$. Then  the concentration estimates of Theorem~\ref{T:CorrLinear} give that
$$
	\lim_{K\to\infty}  \max_{Q\in \Phi_K}|L(f,Q)- f(Q)\cdot Q^{it}\cdot B(f,K)|=0,
$$
where
$$
B(f,K):=\exp(2F_{2N}(f,K))\cdot \overline{ \exp(F_{N}(f,K))}\cdot
\lim_{N\to\infty}\E_{m,n\in [N]} \,  (m^2-n^2)^{it} \cdot m^{-it}\cdot \overline{f(2n)}.
$$
Although the expressions $B(f,K)$ are much more complicated than the expressions $A(f)$ in \eqref{E:Af},  this is of little concern to us, what matters is that for $Q\in \Phi_K$ they are constant in $Q$. This allows us to show that as long as $f\neq n^{-it}$ we have
$$
 \lim_{K\to\infty} \E_{Q\in \Phi_K} f(Q)\cdot Q^{it}\cdot B(f,K)=0.
$$
 As before,  using the bounded convergence theorem, we conclude  that
$$
\lim_{K\to\infty} \E_{Q\in \Phi_K} \int_{\CM_p} 	L(f,Q) \,  d\sigma(f)= \int_{\CA} 	B(f) \,  d\sigma(f),
$$
where
$$
\CA:=\{(n^{it})_{n\in\N}\colon t\in \R\}
$$
and
$$
B(f):= \lim_{N\to\infty}\E_{m,n\in [N]} \,  (m^2-n^2)^{it} \cdot (2mn)^{-it}
$$
when $f=n^{it}$.
  To complete the proof, it suffices to prove the positivity of
 $ \int_{\CA} 	B(f) \, d\sigma(f)$. We do this in two different ways below, the second being the more flexible.
\medskip

(Endgame I: Miraculous positivity.) The limit defining $B(f)$ is easily shown to be
$$
\int_0^1 \int_0^1 (x-y)^{it}\cdot (2xy)^{it}\, dx\, dy.
$$
It can be  shown that this double integral has  positive real part for every $t\in \R$ by explicitly evaluating  the integral (we did this computation using some  computer software). Since $\sigma(\{1\})>0$ this implies that  $ \int_{\CA} 	B(f) \,  d\sigma(f)>0$ and completes the proof.
\medskip

However, this miraculous positivity property fails when we replace $2mn$ with $mn$ or $3mn$ and probably only works for  $kmn$  when  $k=2$, i.e.  only for the case of Pythagorean pairs! Next, we will give an alternative, more complicated, but robust method that works in much higher generality (this is needed for equations other than the Pythagorean).
\medskip

(Endgame II: A more robust method.) We introduce some weights that effectively allow us to replace all non-trivial Archimedean characters with $1$. 	For every $\delta\in (0,1/2)$  we let
\begin{equation}\label{E:wd}
	w_\delta(m,n):={\bf 1}_{I_\delta}\big( (m^2-n^2)^i\cdot (mn)^{-i}\big) \cdot {\bf 1}_{m> n}, \quad m,n\in\N,
\end{equation}
where $I_\delta$ is an arc in $\S^1$ around $1$ with length $\delta$.
It is not difficult to show that for all $\delta>0$ we have
\begin{equation}\label{E:mud}
\mu_\delta:=\lim_{N\to\infty} \E_{m,n\in[N]}\,  w_\delta(m,n)>0
\end{equation}
and also
$$
\lim_{N\to\infty} \E_{m,n\in[N]} \sup_{Q\in \N} |w_\delta(Qm+1,Qn)- w_\delta(m,n)|=0,
$$
so for our purposes,  the weights $w_\delta(m,n)$ and $w_\delta(Qm+1,Qn)$ can be used interchangeably. It is also not hard to show that for aperiodic multiplicative functions $f$ we have the vanishing property \eqref{E:PythVanish} even if we insert the weights
$w_\delta(m,n)$ in the average.
Using these properties, we can work out  the previous argument, replacing  the  expressions $L(f,Q)$ in \eqref{E:LfQ} by
\begin{equation}\label{E:LfQd}
	L_\delta(f,Q):=\lim_{N\to \infty} \E_{m,n\in[N]}\, w_\delta(m,n)\cdot  f((Qm+1)^2-(Qn)^2)\cdot \overline{f(2(Qm+1) Qn)}.
\end{equation}
In this way, we reduce  matters to showing that for some sufficiently small $\delta>0$ we have
$$
\int_{\CA} 	B_\delta(f) \,  d\sigma(f)>0,
$$
where
$$
 B_\delta(n^{it}):= \lim_{N\to\infty}\E_{m,n\in [N]} \,  w_\delta(m,n)\cdot  (m^2-n^2)^{it} \cdot (2mn)^{-it}
$$
when $f=n^{it}$.
The advantage of this weighted average is that for any fixed $t\in \R$, because of the definition of the weight in \eqref{E:wd}, if $\delta$ is small enough, then $B_\delta(f)$ is close to $1$. This allows us to show with some fiddling,  that if $\delta$ is small enough, then we have the lower bound
$$
 \int_{\CA} 	B_\delta(f) \,  d\sigma(f) \geq \frac{1}{2} \cdot \mu_\delta\cdot  \sigma(\CA)\geq\frac{1}{2} \cdot \mu_\delta\cdot  \sigma(\{1\})>0,
$$
where $\mu_\delta$ is as in \eqref{E:mud}.
This completes the proof in a robust way that allows us to deal with positivity properties of more general pairs such as $(\ell mn,\ell' (m^2-n^2))$, where $\ell, \ell'\in \N$ are arbitrary.
The reader will find more details  in \cite[Section~4]{FKM23}.


 \subsection{Comparison with the strategy from \cite{FH17}}
 In \cite{FH17} we were able to show partition regularity for the pairs
 $(m^2-n^2, m(m+2n))$ but not for the pairs
   $(m^2-n^2, 2mn)$. Let us briefly explain why. The main technical
   tool used in \cite{FH17} was a structural result for general multiplicative functions which roughly states the following (see  \cite[Theorem~2.1]{FH17} for the detailed statement): For every  $f\in \CM$, $\varepsilon>0$, and $ s\in \N$,   there exist $Q\in \N$ and $C>0$,  such that for every $N\in \N$  we have a decomposition
   $$
   f(n)=f_{st}(n)+f_{un}(n) +f_{er}(n), \quad n\in[N],
   $$
   such that
 \begin{enumerate}
\item \label{I:1}     $|f_{st}(n+Q)-f_{st}(n)|\leq C/N$ for every $n\in [N-Q]$;

 \item \label{I:2}  $\norm{f_{un}}_{U^s[N]}$ is ``very small'' in a sense that will not be made precise here;

  \item  \label{I:3}  $\E_{n\in[N]}\,  |f_{er}(n)|\leq \varepsilon$.
   \end{enumerate}
   In fact, the applications to partition regularity require some
uniformity in the choice of $Q$ and $C$ with respect to $f\in \CM$, but that will not concern us here.

Let us see how we can use these properties to get some positiveness for fixed $f\in \CM$ and averages of the form
   $$
   \E_{m,n\in [N]} \, f((m^2-(Qn)^2)\cdot \overline{f(m(m+2Qn))}.
$$
  Using  properties~\eqref{I:2} and \eqref{I:3}  it is not hard to see  that   these averages
are $O(\varepsilon)$-close to the averages
  $$
\E_{m,n\in [N]}\,  f_{st}(m-Qn)\cdot f_{st}(m+Qn)
\cdot \overline{f_{st}(m)} \cdot \overline{f_{st}(m+2Qn)}.
$$
Assuming for simplicity that $f_{st}(n+Q)=f_{st}(n)$ for every $n\in \N$  the last average is equal to
 $$
\E_{m,n\in [N]} \, |f_{st}(m)|^4.
$$
This positive lower bound is crucial for the argument in \cite{FH17} and with some  maneuvering can be used to prove partition regularity of the patterns $(m^2-n^2, m(m+2n))$.
However, it  does not appear to be possible to  modify this argument  in a simple way in order to deal with the pairs $(m^2-n^2, 2mn)$, the reason being that  if  $m=0$, or $n=0$, the two coordinates of this pair do not coincide.

Instead, one can try to restrict both coordinates $m$ and $n$ to suitable progressions (just as we did in Section~\ref{SS:pair1}) and work with the averages
\begin{equation}\label{E:QQ}
\E_{m,n\in [N]} \, f((Qm+1)^2-(Qn)^2)\cdot \overline{f((Qm+1)(Qn))}.
\end{equation}
Unfortunately, the previous decomposition result seems inadequate for analyzing these expressions because we have no useful control over the error term $f_{\er}$ along the progression $Q\, \N$.
The important advantage of the concentration estimates described in Section~\ref{S:concentration} is that for pretentious multiplicative functions we obtain an error term that is small even when we restrict both variables $m,n$ to arithmetic progressions. This  is crucial
for the purposes of analyzing the asymptotic behavior of the averages \eqref{E:QQ}.
Moreover, we will see in the next subsection that  the concentration estimates described in Section~\ref{SS:concquad} allow us to effectively analyze averages of the form
$$
\E_{m,n\in [N]}\,  f((Qm+1)^2+(Qn)^2)\cdot \overline{f((Qm+1)(Qn))}
$$
for pretentious multiplicative functions. Such a task
was beyond the limits of the methodology used in the article \cite{FH17}.

\subsection{The pair $(2mn, m^2+n^2)$}\label{SS:pair2}
To establish Theorem~\ref{T:MainMulti1} in this case ($\alpha=\beta=1$, $\gamma=2$), and thus deduce partition and density regularity of $x^2+y^2=z^2$ with respect to the variables  $x,z$,  we follow the argument of the previous subsection with    some necessary and rather laborious changes. We introduce the averages
\begin{equation}\label{E:LfQd'}
	L'_{\delta,N}(f,Q):= \E_{m,n\in[N]}\, \tilde{w}_\delta(m,n)\cdot  f((Qm+1)^2+(Qn)^2)\cdot \overline{f(2(Qm+1) Qn)},
\end{equation}
where the weights $\tilde{w}_\delta$ are defined similarly to \eqref{E:wd}, but with $m^2-n^2$ replaced by $m^2+n^2$.
In order to study their asymptotic behavior, as $N\to\infty$, we first show  using  Theorem~\ref{T:CorrSumSquares}   that in the aperiodic case they vanish for every $Q\in \N$.

We   then
 study their asymptotic behavior  when $f$ is pretentious. To simplify their form, we  use the non-linear concentration estimates of Theorem~\ref{T:concquadratic}.
 This is an important difference, since the non-linear concentration estimates were not previously known and turned out to be  much harder to obtain than the linear ones.
  Another technical annoyance is that, contrary to what happened in  the previous subsection, the limit $\lim_{N\to\infty} L'_{\delta,N}(f,Q)$ may not exist, the reason being that the oscillatory factors resulting from the concentration estimates for the quadratic and linear terms do not conveniently cancel.  In any case, we can show that if $f\sim \chi\cdot n^{it}$ for some Dirichlet character $\chi$ and $t\in \R$, then
 \begin{equation}\label{E:Lfa'}
 	\lim_{K\to\infty}  \limsup_{N\to\infty} \max_{Q\in \Phi_K}|L_{\delta, N}'(f,Q)- f(Q)\cdot Q^{it}\cdot B_{\delta, N}(f,K)|=0,
 \end{equation}
where
	$$
B_{\delta, N}(f,K):=  \exp(G_N(f,K))\cdot \overline{\exp(F_N(f,K))}\cdot  \E_{m,n\in [N]}\, \tilde{w}_\delta(m,n)\cdot  (m^2+n^2)^{it}\cdot  m^{-it}\cdot \overline{f(2n)},
$$
and $F_N(f,K)$ and $G_N(f,K)$ are given by  \eqref{E:FNQ} and \eqref{E:FNQsquares}  respectively.
Again, the complexity of the expression $B_{\delta, N}(f,K)$ will not bother us, the only important fact is that for $Q\in \Phi_K$ it is constant in $Q$. This allows us to show that as long as $f\neq n^{-it}$ we have
$$
\lim_{K\to\infty} \limsup_{N\to\infty}|\E_{Q\in \Phi_K} f(Q)\cdot Q^{it}\cdot B_{\delta, N}(f,K)|=0.
$$
After this point, arguing as in the previous subsection,
we get that for some $\delta_0>0$ we have the following positivity property
\begin{equation}\label{E:positiveQ}
\liminf_{K\to \infty} \liminf_{N\to\infty}
\E_{Q\in \Phi_K} \int_{\CM} L'_{\delta_0,N}(f,Q)\,  d\sigma(f) >0,
\end{equation}
 where $L'_{\delta_0,N}(f,Q)$ is given by \eqref{E:LfQd'}. It is slightly annoying that $\E_{Q\in \Phi_K}$ and
 and $\liminf_{N\to\infty}$ are not interchangeable (in the previous subsection, the existence of the limit as $N\to\infty$ allowed us to do this). To overcome this and finish the argument, we argue as follows. We note that using \eqref{E:positiveQ}
and the fact that $\tilde{w}_\delta(m,n)$ is approximately equal to
$\tilde{w}_\delta(Qm+1,Qn)$, we get
 that there exist $K_0\in \N$ and $Q_N\in \Phi_{K_0}$, $N\in\N$, such that
\begin{equation}\label{E:positiveQN}
 \liminf_{N\to\infty}
 \E_{m,n\in [N]}\,  I(Q_Nm+1,Q_Nn)=:a>0,
\end{equation}
where
$$
I(m,n):=\int_{\CM} \tilde{w}_\delta(m,n)\cdot  f(m^2+n^2)\cdot \overline{f(2m n)} \, d\sigma(f).
$$
Since $I(m,n)\geq 0$ for every $m,n\in\N$ and $Q_N$ belongs to a finite set bounded by, say, $M$, we deduce from \eqref{E:positiveQN} that
$$
\liminf_{N\to\infty}
\E_{m,n\in [N]}\,  I(m,n)=a/M^2>0.
$$

The reader will find the details  of this argument in \cite[Section~6]{FKM23}.

 \subsection{The pair $(mn, m^2+dn^2)$}
The two main new ingredients needed to cover this more general  case are Theorem~\ref{T:CorrGeneral}, which establishes the vanishing property of correlations in the case of aperiodic multiplicative functions, and Theorem~\ref{T:concquadraticgeneral}, which establishes concentration estimates for pretentious multiplicative functions for the binary quadratic form $m^2+dn^2$ when it is irreducible. Modulo these two results, which require considerable additional effort to prove, the required positivity property goes along the lines described in the previous subsections. The details appear in \cite[Section~5]{FKM24}.

We also remark that using a similar approach, we can cover all   pairs $(P_1(m,n),P_2(m,n))$ where at
least one of the $P_1,P_2$ is a reducible binary quadratic form that is not a square.

\subsection{The pair $(m^2+2n^2, m^2-2n^2)$}   We need to cover this pair to prove partition and density regularity of
$x^2+2y^2=z^2$ with respect to the variables $x,z$,
(see the parametrization in Section~\ref{SS:parametric}).
The situation is similar for the equation $x^2+y^2=2z^2$ with respect to any pair of variables.
In order to prove the positivity property of Theorem~\ref{T:MainMulti2} for $\alpha=\alpha'=0$, $\beta=2$, $\beta'=-2$, we need to study  the correlations
\begin{equation}\label{E:fpm2}
\E_{m,n\in [N]} \, f(m^2+2n^2)\cdot \overline{f(m^2-2n^2)}
\end{equation}
for arbitrary $f\in \CM$.  It is expected that asymptotic properties of $f$ along $m^2+ 2n^2$ and $m-2n^2$ are only affected by its values on primes in $\CP_1$ and $\CP_{2}$ respectively, where
\begin{equation}\label{E:P12}
\CP_1:=\{p\equiv 1, 3\!\pmod{8}\}, \qquad \CP_{2}:=\{ p\equiv 1, 7\!\!\pmod{8}\}.
\end{equation}
 (Using the notation in \eqref{E:Pd}, we should have denoted the two sets by $\CP_2$ and $\CP_{-2}$ respectively; we do not do this to simplify the notation later.)

Analogously to \eqref{D:Pretentious},  if $\CP$ is a subset of the primes  we  let
$$
\mathbb{D}^2_\CP(f,g):=\sum_{p\in \CP} \frac{1}{p} \big(1-\Re(f(p)\cdot \overline{g(p)})\big),
$$
and write $f\sim_{\CP} g$ if $\mathbb{D}_\CP(f,g)<+\infty$.
We  say that $f$ is {\em $\CP$-pretentious}  if $\mathbb{D}_\CP(f,\chi\cdot n^{it})<+\infty$ for some Dirichlet character $\chi$ and $t\in \R$, and   {\em $\CP$-aperiodic} if it is not
$\CP$-pretentious.   Note that for $\CP=\P$ we get the same definition
as in Section~\ref{SS:pretentious}.

\subsubsection{The pretentious case}  Let $\CP_1,\CP_2$ be as in \eqref{E:P12} and  suppose  that $f$ is $\CP_1$-pretentious and $\CP_2$-pretentious. Furthermore, for convenience, suppose that
\begin{equation}\label{E:fpret}
	f\sim_{\CP_1}\chi_1 \quad  \text{and} \quad
		f\sim_{\CP_2}\chi_2
	\end{equation}
	for some Dirichlet characters $\chi_1,\chi_2$. We will ignore Archimedean characters in this discussion, since they can be covered by introducing weights and using a straightforward variant of the argument given in Section~\ref{SS:pair2}.

Our aim is to use again the strategy described in the previous subsections. However,  substantial new complications arise
since we no longer have a reducible polynomial like $mn$ to work with and
this  makes the
$Q$-trick harder  to implement. Moreover,    Theorem~\ref{T:concquadraticgeneral} gives that there exist
multiplicative functions which concentrate on different oscillatory factors along $m^2+2n^2$ and $m^2-2n^2$
that do not  cancel each other. This   casts doubt on whether it is even possible in this case to establish positivity via
concentration estimates.  Fortunately, we will be able to do this  by adapting the $Q$-trick, this time averaging over appropriate ``translates'' rather than ``dilates'' of an appropriate choice of lattice.   Our goal is to ensure that for those multiplicative functions that concentrate
along $m^2+2n^2$ and $m^2-2n^2$
around different oscillatory factors, a certain multiplicative average over $Q$ will vanish, and for the remaining multiplicative functions we can guarantee that they concentrate on the same oscillatory factor, which leads to a favorable cancellation.

Let us explain in  more detail how the previous plan is implemented.
For a given $K\in \N$  and $l_p,l'_p\in [1,3K/2]$, let
\begin{equation}\label{E:Q12}
Q_1:=\prod_{p\equiv 3  \! \!  \!   \! \!  \pmod{8}, \, p\leq K, \, p\nmid q}p^{l_p}, \quad
Q_2:=\prod_{p\equiv 7  \!   \!  \!  \!    \pmod{8}, \, p\leq K, \, p\nmid q}p^{l'_p}.
\end{equation}
Let $q$ be a common period of the mutliplicative functions $\chi_1$ and $\chi_2$. We plan to apply concentration estimates on the lattice
$$
\{(Q_Km+a,Q_Kn+b)\colon m, n\in \N\},
$$
where
$$
Q_K:=\prod_{p\leq K} p^{2K}
$$
and  $a,b\in [Q]$ depend on $Q_1,Q_2$ and satisfy the following two properties
\begin{equation}\label{E:ab1}
	a^2+2b^2\equiv a^2-2b^2\equiv 1\pmod{q}
\end{equation}
and
\begin{equation}\label{E:ab2}
	(a^2+2b^2,Q)=Q_1, \quad (a^2-2b^2,Q)=Q_2.
\end{equation}
To show that a choice of $a,b$ that satisfy \eqref{E:ab1} and \eqref{E:ab2}  is possible, we
use the Chinese remainder theorem and the fact that for $p\equiv 3\pmod{8}$ the
congruences $p^{l_p} \mid \mid a^2+2b^2$  and $p\nmid a^2-2b^2$ have a solution, and similarly,  for $p\equiv 7 \pmod{8}$ the
congruences $p^{l_p} \mid \mid a^2-2b^2$  and $p\nmid a^2+2b^2$ have a solution.

 Suppose now that $K$ is sufficiently large, $Q_1, Q_2$ are as in \eqref{E:Q12}, and  $a=a_{Q_1,Q_2},b=b_{Q_1,Q_2}$ satisfy  \eqref{E:ab1} and \eqref{E:ab2}.  For convenience we let
 $$
 P_1(m,n):=m^2+2n^2,\quad P_2(m,n):=m^2-2n^2.
 $$
 Using  a variant of Theorem~\ref{T:concquadraticgeneral} that allows us to work with
 lattices of the form $Q\cdot \N\times \N+ (a,b)$ we get
  that for $j=1,2$ and  $Q=Q_K$ the values of
  $f(P_j(Qm+a,Qn+b))$ are concentrated at
\begin{equation}\label{E:fQj}
f(Q_{j})\cdot \chi_j(P_j(a,b)/Q_j) \cdot \exp(G_{j,N}(f,K)),
\end{equation}
 (we used  that   $Q_j=(P_j(a,b), Q)$, $j=1,2$, and  $q\prod_{p\leq K}p$ divides $Q/Q_{j}$)
where
$$
	G_{j,N}(f,K):= 2\, \sum_{ K< p\leq N, \, p\in \CP_j}\, \frac{1}{p} \,(f(p)\cdot \overline{\chi_j(p)}\cdot -1), \quad j=1,2.
$$
Note that for $j=1,2$ we have    $(Q_j,q)=1$, hence  $|\chi_j(Q_j)|=1$, and also    since  \eqref{E:ab1} holds, we have
$\chi_j(P_j(a,b))=1$. It follows that the expression in \eqref{E:fQj} is equal to
$$
f(Q_j)\cdot \overline{\chi_j(Q_{j})} \cdot \exp(G_{j,N}(f,K)).
$$
Hence,  $f(P_1(Qm+a,Qn+b))\cdot \overline{f(P_2(Qm+a,Qn+b))}$ concentrates at
\begin{equation}\label{E:concQ12}
 f(Q_1)\cdot \overline{\chi_1(Q_1)} \cdot \exp(G_{1,N}(f,K))
 \cdot \overline{f(Q_2)} \cdot \chi_2(Q_2) \cdot \overline{\exp(G_{2,N}(f,K))}.
\end{equation}
We consider two cases.

\medskip

{\bf Case 1. } Suppose that  for  primes $p$ that satisfy $ p\nmid q$ we have
$$
	f(p)\neq \chi_1(p) \text{ for some  } p\equiv 3 \! \! \! \pmod{8} \quad \text{or }
	\quad f(p)\neq \chi_2(p) \text{  for some  } p\equiv 7 \! \!\!\pmod{8}.
$$
Then the previous facts easily imply that
\begin{equation}\label{E:LfK}
\lim_{K\to\infty} \limsup_{N\to\infty}|\E_{Q_1\in \Phi_{1,K}, Q_2\in \Phi_{2,K}}\,
L_N(f,K,Q_1,Q_2)|=0,
\end{equation}
where
\begin{multline}\label{E:LKQ12}
L_N(f,K,Q_1,Q_2):=\\
\E_{m,n\in [N]}\, f(P_j(Q_Km+a_{Q_1,Q_2},Q_Kn+b_{Q_1,Q_2}))\cdot \overline{f(P_j(Q_Km+a_{Q_1,Q_2},Q_Kn+b_{Q_1,Q_2})}),
\end{multline}
and  for $j=1,2$
		$$
	\Phi_{j,K}:=\Big\{ \prod_{p\in \CP_j\cap [K]} p^{l_p}\colon K< l_p\leq  3K/2\Big\}.
	$$

\medskip

{\bf Case 2.} Suppose that for     primes $p$ that satisfy $p\nmid q$ we have
 \begin{equation}\label{E:fid12}
f(p)= \chi_1(p) \text{ for all  } p\equiv 3 \! \! \! \pmod{8} \quad \text{and }
\quad f(p)= \chi_3(p) \text{  for all  } p\equiv 7 \! \!\!\pmod{8},
\end{equation}
  then the expression in \eqref{E:concQ12} equals
 \begin{equation}\label{E:concQ12'}
 	 \exp(G_{1,N}(f,K))
 	\cdot \overline{ \exp(G_{2,N}(f,K))}.
 \end{equation}
Using
\eqref{E:fid12} we get  that for $K$ large enough  and $N>K$ we have
$$
G_{1,N}(f,K)=G_{2,N}(f,K).
$$
Thus the expression in \eqref{E:concQ12'} equals
$$
 \exp(2\Re(G_{1,N}(f,K))).
$$
Since our standing assumption \eqref{E:fpret} is that $f\sim_{P_1} \chi_1$, this last  expression  approaches $1$  as $N\to\infty$ uniformly on $N>K$.
Combining the above we get
\begin{equation}\label{E:LfK'}
	\lim_{K\to\infty} \limsup_{N\to\infty} \sup_{Q_1\in \Phi_{1,K}, Q_2\in \Phi_{2,K} }
	|L_N(f,K,Q_1,Q_2)-1|=0,
\end{equation}
where $L_N(f,K,Q_1,Q_2)$ is as in \eqref{E:LKQ12}.

  Combining \eqref{E:LfK} and  \eqref{E:LfK'} we get the  positiveness properties needed to prove Theorem~\ref{T:MainMulti2}.

\subsubsection{The aperiodic case}\label{SS:apChowla}
 In the case where the multiplicative function $f\in \CM$ is $\CP_1$-aperiodic or $\CP_2$-aperiodic, we expect the averages \eqref{E:fpm2} to vanish as $N\to \infty$. Unfortunately,   this property  turns out to be very  difficult to establish because, unlike the cases treated in the previous subsections, we have no way of bounding these averages by a usable substitute of the Gowers uniformity seminorms of
 the multiplicative function $f$.  In fact, the needed vanishing property is not even known for the Liouville function, i.e. it is not known that
	$$
	\lim_{N\to \infty}\E_{m,n\in[N]}\, \lambda(m^2+2n^2)\cdot
	\lambda(m^2-2n^2)=0.
	$$
For the applications we have in mind we  would be equally happy to show this for logarithmic averages in place of Ces\`aro.

\begin{definition}\label{D:CorVan}
	If   $P\in \Z[m,n]$ is an  irreducible  binary quadratic form with  discriminant $d_P$, we let
	$$
		\CP_P:=\Big\{p\in \P\colon  \legendre{d_P}{p}=1 \Big\}.	
	$$
	We say that the irreducible binary quadratic forms   $P_1,P_2$ are  {\em good for vanishing of correlations of aperiodic multiplicative functions,} whenever   the following property holds: If  $f_1,f_2\colon \N\to \S^1$ are completely multiplicative functions  such that either  $f_1$ is $\CP_{P_1}$-aperiodic or $f_2$ is $\CP_{P_2}$-aperiodic,
	then 	for every $Q\in \N$ and $a,b\in \Z_+$ we have
	$$
	\lim_{N\to\infty}\E_{m,n\in [N]} \,{\bf 1}_S(Qm+a,Qn+b)\cdot  f_1(P_1(Qm+a,Qn+b)) \cdot f_2(P_2(Qm+a,Qn+b))=0,
	$$
	where $S:=\{m,n\in\N\colon P_1(m,n)>0,P_2(m,n)>0\}$.
\end{definition}

	\begin{conjecture}\label{C:2quadratic}	
	If $P_1,P_2\in \Z[m,n]$ are irreducible binary quadratic forms that are not multiples of each other,
	then they are good for vanishing of correlations of aperiodic multiplicative  functions.
\end{conjecture}
  If this conjecture is true, we show in  \cite[Section~6]{FKM24} how Theorem~\ref{T:MainMulti2} follows using the ideas summarized in the previous subsection.

\section{More challenges}
In this section we present some open questions linked to the topics already discussed.
\subsection{Vanishing of correlations along  two irreducible quadratics}
The following problem covers   a special case of Conjecture~\ref{C:2quadratic} that seems very challenging:
\begin{problem}
	If $\lambda$ is the Liouville function, show that
	$$
	\lim_{N\to \infty}\E_{m,n\in[N]}\, \lambda(m^2+n^2)\cdot
	\lambda(m^2+2n^2)=0.
	$$
\end{problem}
For the applications of  this article we could replace the Ces\`aro averages
$\E_{m,n\in[N]}$ with the logarithmic averages $\frac{1}{(\log{N})^2}\sum_{m,n\in[N]} \frac{1}{mn}$.
In fact, it appears that  no explicit example of two irreducible binary quadratic forms
$P_1,P_2$ is known  for which
 $$
 \lim_{N\to \infty}\E_{m,n\in[N]}\, \lambda(P_1(m,n))\cdot
 \lambda(P_2(m,n))=0,
 $$
or a similar property with logarithmic averages.

\subsection{Partition regularity of ``most'' generalized Pythagorean triples }
Conjecture~\ref{C:2quadratic} asks to show that for all Rado triples $(a,b,c)$ the equation \eqref{E:pythgen} is partition regular. At the moment this seems very hard
to prove, so even showing that we have partition regularity for a positive proportion, or even better, for most Rado triples, would be desirable.
\begin{problem}\label{P:mostabc}
	Show that for ``most'' choices of Rado triples $(a,b,c)\in \N^3$ the equation
	$$
	ax^2+by^2=cz^2
	$$
	is partition regular.
\end{problem}
We do not make explicit what ``most'' means, but let us say that a subset of Rado's triples with
relative density $1$ with respect to the standard density, would qualify. At the moment it is not even known that there is a single Rado triple for which we have partition regularity.

\subsection{Multiple recurrence results for multiplicative actions} \label{SS:multi}
Although the methods of this article are particularly well suited to partition regularity problems of pairs,
they seem inadequate on their own to deal with problems of triples.
The reason is that we lack a representation result for sequences such as
$$
(x,y,z)\mapsto d(x^{-1}\Lambda\cap y^{-1}\Lambda\cap  z^{-1}\Lambda),
$$
which was our starting point in our approach to partition regularity of pairs.  It is understood that a useful representation
similar to the Bochner-Herglotz theorem should be very hard to prove, and that new sequences, not just multiplicative functions, should be
used.


An alternative approach, which has been very successful for translation-invariant patterns of arbitrary length, is to use the formalism of ergodic theory to provide
sufficient statements for partition regularity and then to develop the necessary toolbox in ergodic theory to prove them. The sufficient statements in ergodic theory are described next.

We say that $(X,\mu,T_n)$ is a {\em multiplicative action}, if
$(X,\mathcal{X},\mu)$ is a probability space and  for $n\in \Z$ the transformations $T_n\colon X\to X$ are invertible, measure-preserving, and satisfy the relation $T_0=T_1:=\text{id}$ and $T_m\circ T_n=T_{mn}$ for all $m,n\in \Z$.
Using a variant of the correspondence principle of Furstenberg~\cite{Fust}  that appears in  \cite{Be05}, it turns out that in order to show that the equation $x^2+y^2=z^2$ is partition regular, it suffices to prove the following multiple recurrence result:
\begin{problem}\label{P:Pythagoras}
	Let $(X,\mu,T_n)$ be a multiplicative action an $A$ a measurable set such that  $X=\bigcup_{j=1}^kT_j^{-1}A$  for some $k\in \N$. Show that
\begin{equation}\label{E:Aj}
\liminf_{N\to\infty}\E_{m,n\in [N]}\, \mu(T^{-1}_{m^2-n^2}A \cap T^{-1}_{2mn}A\cap T^{-1}_{m^2+n^2}A)>0.
\end{equation}
\end{problem}
Unlike the partition regularity results for pairs, such as those in Theorem~\ref{T:PairsDensityParametric}, we cannot expect to have
density regularity results for the equation
$ax^2+by^2=cz^2$ when $a+b\neq c$.
If, say, $a=b=c=1$, an example by Bergelson shows that for suitably chosen
chosen $\alpha\in \R$, the set
$\Lambda:=\{n\in \N\colon \{n^2\alpha \}\in [1/5,2/5)\}$ has positive multiplicative density but does not contain a Pythagorean triple (this problem already appears for the equation $x+y=z$).
This is the reason why in Problem~\ref{P:Pythagoras} we cannot expect to have a multiple recurrence property for  any set $A$ with positive measure.
On the other hand, if the action is finitely generated, i.e. the set $\{ T_p,p\in \P\}$ is finite, we expect that \eqref{E:Aj} holds for  all sets $A$ with positive measure.


In order to show that the equation $x^2+y^2=2z^2$ is partition regular (in fact, density regular with respect to multiplicative density), it suffices to prove the following:
\begin{problem}
	Let $(X,\mu,T_n)$ be a   multiplicative action and let $A\in \mathcal{X}$ with $\mu(A)>0$. Show that
$$
\liminf_{N\to\infty}\E_{m,n\in [N]}\, \mu(T^{-1}_{m^2-n^2+2mn}A \cap T^{-1}_{m^2-n^2-2mn}A \cap T^{-1}_{m^2+n^2}A)>0.
$$	
\end{problem}

The problem is open even if the action is finitely generated and we  omit one of the three terms in the intersection.

\subsection{Higher degree homogeneous polynomials}\label{Q:higherdegree}
As we have seen in Section~\ref{SS:conentrationgeneral}, for every irreducible binary quadratic form $P(m,n)$ all pretentious multiplicative functions enjoy strong concentration properties. Does a similar fact hold for not necessarily homogeneous polynomials or higher degree homogeneous polynomials?  By slightly modifying an example of Klurman  \cite[Lemma~2.1]{Kl17}  one can construct pretentious (and aperiodic) multiplicative functions
$f\colon \N\to \{-1,1\}$ so that the averages
$$
\E_{n\in [N]} \, f(n^2+1)
$$
do not converge.  It follows that in this case we cannot have concentration estimates
similar to those given in Theorem~\ref{T:fv}. Although,  as we have shown in Theorem~\ref{T:concquadratic},  such counterexamples do not work  for, say, $P(m,n)=m^2+n^2$,
it seems likely that they can be found for higher degree irreducible homogeneous polynomials, but this is not so easy to achieve.

\begin{problem}
	Let $a\in \Z$ and $d\geq 3$ be such that the polynomial $m^d+an^d$ is irreducible.
	Show that  there exists a pretentious completely multiplicative function $f\colon \N\to \{-1,1\}$ such that the averages
	$$
	\E_{m,n\in [N]} \, f(m^d+an^d)
	$$
	do not converge as $N\to \infty$.
\end{problem}

We say that a sequence $a\colon \N^2\to \Z$ is {\em good for single recurrence of multiplicative actions} if the following holds: For every multiplicative action $(X,\mu,T_n)$ and $A\in \CX$ with $\mu(A)>0$, we have
$$
\mu(A\cap T^{-1}_{a(m,n)}A)>0
$$
for some $m,n\in \N$ such that $a(m,n)\neq 0$.
Is  the polynomial $P(m,n)=m^3+2n^3$ good for single recurrence of multiplicative actions? The following is a more general problem:
\begin{problem}
Let $P\in \Z[m,n]$ be a homogeneous polynomial  with $P(1,0)=1$.  Is it true that  for every 	 multiplicative action  $(X,\mu,T_n)$ and $A\in \mathcal{X}$ with $\mu(A)>0$, we have
	$$
\liminf_{N\to\infty}\E_{m,n\in [N]}\,	\mu(A \cap T^{-1}_{P(m,n)}A)>0 \, ?
	$$	
\end{problem}
We cannot remove the assumption $P(1,0)=1$ since the set  $\{2n^2\colon n\in \N\}$ can easily be seen to be bad for measurable recurrence \cite[Example~3.11]{DLMS23}.


It follows from \cite[Theorem~2.1 and Example~5.15]{DLMS23} that the answer is positive when  $\deg(P)\leq 2$.
It seems very likely that, using the methodology  developed in \cite{FKM24}, we can also get a positive  answer  for  reducible polynomials of degree $3$. But for  $P(m,n)=m^3+2n^3$,  it is not even clear how to show that for  every completely multiplicative function  $f\colon \N\to \{-1, 1\}$ we have $f(m^3+2n^3)=1$ for a set of $m,n\in \N$ with positive lower  density.

		\end{document}